\input amstex
\documentstyle{amsppt}
\input epsf
\nologo
\nopagenumbers
\NoBlackBoxes
\hcorrection{1.5cm}
\topmatter
\title Specialization of polynomial covers of prime degree\endtitle
\author Leonardo Zapponi\endauthor
\date October 2002\enddate
\abstract Let $K$ be a complete field of unequal characteristics $(0,p)$ . The aim of this paper is to describe the semi-stable models for covers $\bold P^1_K@>>>\bold P^1_K$ of degree $p$, unramified outside $r\leq p$ points and totally ramified above one of them, under the assumtion that the ramification locus has a particular reduction type (which always occurs if $r\leq 4$). We are principally concerned with the minimal semi-stable models which separate the ramified fibers.
\endabstract
\leftheadtext{Leonardo Zapponi}
\rightheadtext{Specialization of polynomial covers}
\address Ecole Polytechnique F\'ed\'erale de Lausanne\endaddress
\email leonardo.zapponi\@epfl.ch\endemail
\urladdr http://alg-geo.epfl.ch/\~{}marina/leonardo.html\endurladdr
\endtopmatter
\document

\specialhead \S0 Introduction\endspecialhead

The study of the stable models for covers between algebraic curves over a complete discrete valuation ring is a subject which has been intensively developped during these last years. Roughly speaking, the problem can be summarized as follows: starting from a ramified cover $\beta:C@>>>D$ between non-singular projective curves defined over a complete discrete valuated field $(K,v)$ of unequal characteristics, one would like to construct a semi-stable model $\Cal C@>>>\Cal D$ of this cover over its valuation ring $R$. The general theory, and in particular the semi-stable reduction theorem, asserts that such a model always exists, up to a finite extension of the base field (cf. [BLR] and [R1]). Its uniqueness is not ensured, but  if we require that the model separates the ramified locus (that is, the ramification points specialize to pairwise distinct points), then there exists a minimal semi-stable model for the cover, which is unique, up to isomorphism. Moreover, this model behaves well under (finite) base change, so that it can be considered as an important birational invariant attached to the cover. In particular, there is a well defined notion of reduction type: the number of irreducible components of the special fiber, their intersection graph and the thicknesses of the singular points (cf. section 1).  

Historically, the investigation started with the study of tame covers, i.e. covers such that the residual characteristic $p$ of $R$ does not divide the order of their monodromy group $G$. This first approach is somehow more simple, and one can use the results of [G], for example. The next case is where $p$ divides $|G|$. There are no complete results on this direction. Raynaud [R2] treates the case where $G$ has a $p$-Sylow subgroup of order $p$ and the cover is unramified outside three points. These results are sharpened by Wewers in [W]. Finally, Lehr [L] and Sa\"\i di [S], study the stable model for a $p$-cyclic cover in full generality. In order to have a more complete reference on the subject, see also the results of Green-Matignon and Henrio.

In this paper, we concentrate on the covers $\bold P^1@>>>\bold P^1$ of prime degree $p$. The only restrictions we make on the cover are: there is at least one totally ramified point (wild ramification) and the branch locus has a particular reduction type (cf.~section 5). Note that these covers need not to be Galois. It turns out that the stable models only depend on the ramification data and on the $v$-adic distance between the ramified points. There are  four cases: if the branch locus has good reduction, we essentially recover Wewers' results. The real interesting case is where the branch points have bad reduction. We choose two branch points $P$ and $Q$ and consider the semi-stable model $\Cal C@>>>\Cal D$ of the cover which separates these points. There are three cases: $P$ and $Q$ are $v$-adically far from each other, $P$ and $Q$  have medium distance from each other (the critical case) and $P$ and $Q$ are $v$-adically close to each other. In the first  case, $P$ and $Q$ will specialize into two different irreducible components $C_1$ and $C_2$, and there will be only onle irreducible component lying above each $C_i$. In the second case, $P$ and $Q$ will specialize in the same irreducible component $C$ and there will be only one component above $C$. In the last case, $P$ and $Q$ will again specialize into the same component $C$, but many irreducible components will lye over it. 

The strategy is to start with a natural model for the cover $\bold P^1@>>>\bold P^1$ over $R$ and then to blow it up several  times untill the ramification points are separated. This construction can be carried out explicitely because the covering curve is the projective line. This techniques can be carried over to the case where the degree of the cover is composite and even to the case where the curve $C$ has (potentially) good reduction, but we do not do it in this paper.

\specialhead \S1 Some notation and basic facts\endspecialhead

Throughout this paper, $(K,v)$ will denote a complete valuated field of unequal characteristics $(0,p)$. Its valuation ring $R$ is a complete discrete valuation ring with maximal ideal $\frak p$ and residue field $k=\Cal O_K/\frak p$. The valuation $v:K^*@>>>\bold Q$ will be normalized by the condition $v(p)=1$, so that it will extend the usual $p$-adic valuation. In particular $v(K^*)=e^{-1}\bold Z$, where the positive integer $e$ is the absolute ramification index of $K$, that is $pR=\frak p^e$.

Let $C$ be a projective, irreducible, non-singular curve of genus $g$ defined over $K$ and consider a finite subset $S$ of $C(K)$ of cardinality $n\geq 0$ (for $n=0$, we will set $S=\emptyset$). The pair $(C,S)$ is called a {\it pointed curve of type} $(g,n)$. We will say that $(C,S)$ is {\it hyperbolic} if the inequality $2g-2+n>0$ holds. A {\it semi-stable model} for $(C,S)$ over $R$ is a proper and flat scheme $\Cal C/R$ satisfying the following conditions:
\roster
\item The generic fiber $\Cal C_K=\Cal C\otimes_RK$ is isomorphic to $C$.
\item The special fiber $\Cal C_k=\Cal C\otimes_Rk$ is reduced and has only ordinary double points as singularities.
\item The elements of $S$ specialize to pairwise distinct non-singular points of $\Cal C_k$.
\endroster 

If $S=\emptyset$, then this last condition must be ignored. We will say that $\Cal C$ is {\it stable} if, for any irreducible component $X$ of $\Cal C_k$, we have the relation $2g(X)-2+n(X)+s(X)>0$, 
where $g(X)$ denotes the genus of $X$, $n(X)$ is the cardinality of the subset of $S$ specializing in $X$ and $s(X)$ is the number of singular points of $X$. Semi-stable models behave well under base change, i.e. if $L/K$ is a finite extension, and $\Cal C/R$ is a semi-stable model of $(C,S)$ over $R$, then $\Cal C\otimes_{R}R'$ is a semi-stable model for $(C,S)$ over $R'$, where $R'$ is the ring of integers of $L$. 

The semi-stable reduction theorem asserts that any pointed curve $(C,S)$ over $K$ admits a semi-stable model over the valuation ring  $R'$ of a finite extension $L$ of $K$. Moreover, if $(C,S)$ is hyperbolic, then there exists a stable model $\Cal C/R'$, which is unique, up to $R'$-isomoprhism. The same kind of result holds for covers: more precisely,  given a finite cover $\beta:C@>>>D$ of projective, non-singular curves over $K$ and a finite subset  $S$ of $D(K)$, there exist a finite extension $L$ of $K$, two proper and flat schemes $\Cal C,\Cal D/R'$ and a finite morphism $\Cal C@>>>\Cal D$ such that the following conditions hold:

\roster
\item $\Cal C$ (resp. $\Cal D$) is a semi-stable model for the pointed curve $(C,\beta^{-1}(S))$ (resp. for the pointed curve $(D,S)$) over $R'$.
\item The generic morphism $\Cal C_L@>>>\Cal D_L$ is isomorphic to $\beta:C@>>>D$.
\item The specialized morphism $\Cal C_l@>>>\Cal D_l$ maps non-singular points to non-singular points (here, $l$ denotes the residue field of $L$).
\endroster

 In particular there exists a {\it minimal semi-stable model}  satisfying the above properties, which is unique, up to $R'$-isomorphism. Since  the behaviour of its special fiber is stable  under finite base change, it can be considered as an important birational invariant of the cover. In this paper, we are not concerned with rationality questions, that's why we will always suppose that the semi-stable model is already defined over $R$. 

We will end this section by introducing an important invariant attached to a semi-stable model $\Cal C/R$: the completion of the local ring of $\Cal C$ at a singular point $P$ of its special fiber is isomorphic to the power series ring $R[\![X,Y]\!]/(XY-Z)$, with $Z\in\frak p$. The {\it thickness} of $P$ is the valuation of $Z$, which only depends on $P$. Moreover, our normalization for the valuation $v:K^*@>>>\bold Q$ implies that the thickness of a singular point does not change after a finite base change. Finally, if $\beta:\Cal C@>>>\Cal D$ is a semi-stable model for a cover and $P$ is a singular point of the special fiber $\Cal C_k$ of $\Cal C$ of thickness $\nu$, then $Q=\beta(P)$ will be a singular point of the special fiber $\Cal D_k$ of $\Cal D$, of thickness $e_P\nu$, where $e_P$ is the ramification index at $P$ of the secialized cover $\Cal C_k@>>>\Cal D_k$.

We refer to [BLR] for a complete introduction to the theory of semi-stable curves and their application to the study of algebraic covers of curves. Two fundamental papers on the subject are [R1] and [R2].

\specialhead \S2 Polynomial covers. Normalized models\endspecialhead
 
The main purpose of this paper is to construct semi-stable models for {\it polynomial covers} of prime degree $p$ over $K$ (i.e.  covers $\beta:\bold P^1_K@>>>\bold P^1_K$ of degree $p$ and totally ramified above one point) under certain assumptions on the reduction of the branch locus. We are principally concerned with the minimal models $\Cal X@>>>\Cal Y$ which separates the ramification locus. From now on, we will say that two such covers $\beta_1$ and $\beta_2$ are {\it isomorphic, or equivalent (over $K$)} if there exist two elements $\sigma,\tau\in\text{PGL}_2(K)$ such that $\sigma\circ\beta_1=\beta_2\circ\tau$, i.e. the following diagram
$$\CD\bold P^1_K&@>\,\,\,\tau\,\,\,>>&\bold P^1_K\\
@V\beta_1VV& &@V\beta_2VV\\
\bold P^1_K&@>\,\,\,\sigma\,\,\,>>&\bold P^1_K\endCD$$ 
is commutative. In this case, we easily see that the special fibers of the stable models of the covers $\beta_1$ and $\beta_2$ have the same behaviour. We will start by constructing a ``good'' model of the cover over $K$. Let's first of all reduce to the affine case: since $\beta:\bold P^1_K@>>>\bold P^1_K$ is totally ramified above one point, we obtain a finite morphism $$\beta:\bold A^1_K=\text{Spec}\left(K[X]\right)@>>>\bold A^1_K=\text{Spec}\left(K[T]\right)$$
In terms of affine algebras, it corresponds to an injection $\beta^*:K[T]@>>>K[X]$, which is uniquely determined by the image $\beta(X)$ of $T$, which is a polynomial (that's why we are speaking about polynomial covers). From now on, by eventually enlarging $K$, we will assume that the ramification locus of $\beta$ is contained in $\bold P^1(K)$.
 For any $\lambda\in K$, we have a unique factorization (in $\overline K$)
$$\beta(X)-\lambda=c\prod_{x\in\beta^{-1}(\lambda)}(X-x)^{e_x}$$
where $c\in K$ does not depend on $\lambda$ and $e_x$ is the multiplicity of $X-x$ in the factorization of $\beta(X)-\lambda$. We clearly have 
$$\sum_{x\in\beta^{-1}(\lambda)}e_x=p$$
Furthermore,  the Riemann-Hurwitz formula gives $$\sum_{\lambda\in S}n_\lambda=(r-2)p+1$$
where $r$ is the cardinality of the branch locus $B=S\cup\{\infty\}$ of the cover and $n_\lambda$ is the cardinality of the fiber $\beta^{-1}(\lambda)$. Up to equivalence, and after a finite extension of $K$, we can reduce to the following case:

\roster
\item The finite branch locus $S=\{\lambda_1,\dots,\lambda_{r-1}\}=\beta\left(\{x\,\,\in K\,\,|\,\,\beta'(x)=0\}\right)$ is contained in $R$ and $\{0,1\}\subset S$.
\item  The polynomial $\beta(X)$ is monic (i.e. $c=1$)  and $\beta(0)=0$.
\endroster

\noindent A polynomial $\beta(X)\in K[X]$ satisfying the above conditions will be called is a {\it normalized model} for the cover. One easily checks that there exist finitely many normalized models associated to the same polynomial cover. More generally, we will say that $\beta(X)$ is a {\it semi-normalized} model if it satisfies the condition (2) above and the following property (which is weaker than the condition (1)):

\roster
\item[3]  The finite branch locus $S$ is contained in $R$ and there exists $\lambda\in R^*$ such that $\{0,\lambda\}\subset S$.
\endroster

Remark that there exists infinitely many semi-normalized models associated to the same cover.

\specialhead \S3 Construction of the fundamental model\endspecialhead

We will now start the construction of the semi-stable model associated to a polynomial cover of degree $p$. For any polynomial $h(X)\in R[X]$, we will denote by $\overline h(X)$ its canonical image in $k[X]$. The importance of the (semi-)normalized model introduced in the previous paragraph is related to the following lemma:

\proclaim{3.1 Lemma}\ Let $ h(X)\in K[X]$ be a monic polynomial of degree $p$ such that $h(0)=0$. If the associated cover $h:\bold A^1_K@>>>\bold A^1_K$ is unramified outside a finite set $S\subset R$, then $ h(X)\in R[X]$ and $\overline h(X)=X^p$.\endproclaim

\demo{Proof} The case $ h(X)=X^p$ is immediate. If $ h(X)\neq X^p$ then its Newton polygon is not reduced to a vertical line. Suppose that the last  segment has  positive slope, i.e. that there exists a root $x$ of $ h(X)$ with (minimal) negative valuation. The  degree of $ h(X)$ is  equal to $p$ and $ h(0)=0$, so that  its Newton polygon coincides with the  Newton polygon associated to $X h'(X)$, exept for the last segment (see the following picture). 

\vskip.6cm
 
\epsfysize=4cm

\centerline{\epsfbox{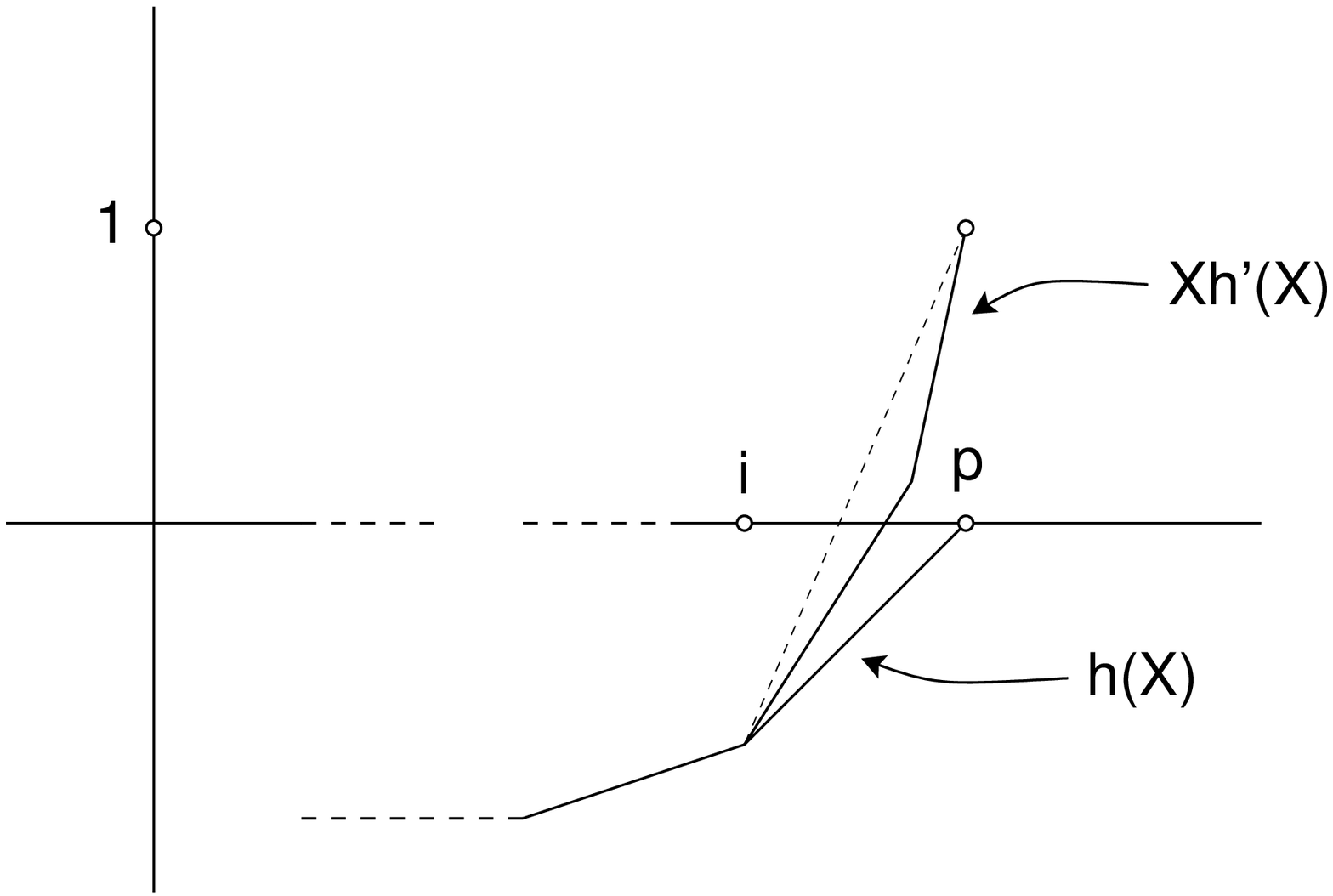}}

\vskip.6cm 

\noindent In particular, we see that there is a root $y$ of $ h'(X)$ (a ramified point) such that its valuation is strictly less than the valuation of $x$ (and thus, of any root of $ h(X)$). We obtain $v( h(y))=pv(y)<0$, which is absurd since  $ h(y)\in S\subset R$. We then have $ h(X)\in R[X]$ and since the (finite) branch points belong to $R$, we see that all the roots of $h'(X)$ belong to $R$. The leading coefficient of $h'(X)$ is equal to $p$, so that $ h'(X)=pg(X)$, with $ g(X)\in R[X]$. This gives $\overline h'(X)=0$, which leads to $\overline h(X)=X^p$, since $ h(X)$ is monic, of degree $p$ and satisfies $ h(0)=0$.\qed \enddemo

 If $\beta(X)\in K[X]$ is a (semi-)normalized model associated to a  polynomial cover of degree $p$ (as defined above), then it satisfies the hypothesis of Lemma 3.1, so that it belongs to $R[X]$ and we have the identity $\overline\beta(X)=X^p$. This will be our starting point for the construction of the semi-stable model associated to the cover. Remark that our assumptions on $K$ imply that, for any $\lambda\in S$, the polynomial $\beta(X)-\lambda\in K[X]$ totally splits. By construction, we moreover have $\beta^{-1}(S)\subset R$. Setting $\Cal C_0=\text{Spec}(R[X])\cong\bold A^1_{R}$ and $\Cal D_0=\text{Spec}(R[T])\cong\bold A^1_{R}$,  we then obtain  a morphism $$\Cal C_0@>>>\Cal D_0$$
 In terms of affine algebras, it corresponds to the injection $R[T]@>>>R[X]$ which sends $T$ to $\beta(X)$. It is a finite morphism, since the polynomial $\beta(X)$ belongs to $R[X]$ and is monic. Let $\gamma(X)\in R[X]$ be the polynomial defined by the relation $\gamma(X)=X^p\beta(X^{-1})$. From a practical point of view, we have
$$\gamma(X)=\prod_{x\in\beta^{-1}(0)}(1-xX)^{e_x}$$
where $e_x=e_x(\beta)$ is the multiplicity of $X-x$ in the factorization of $\beta(X)$. The relation $\beta(0)=0$ implies that the degree of $\gamma(X)$ is strictly less than $p$. Moreover, we have $\gamma(0)=1$ and $\overline\gamma(X)=1$. Set $\frak A=R[T_\infty]$ and $\frak B=\frak A[X_\infty]/(X_\infty^p-\gamma(X_\infty)T_\infty)$. The ring $\frak B$ is a finite $\frak A$-module and thus, we obtain a finite morphism $$\Cal C_1@>>>\Cal D_1$$
 where $\Cal C_1=\text{Spec}(\frak B)$ and $\Cal D_1=\text{Spec}(\frak A)\cong\bold A^1_{R}$. One can easily check that the ring $\frak B$ is canonically isomorphic to $R[X_\infty,\gamma(X_\infty)^{-1}]$. In particular, the generic fiber of $\Cal C_1$ is isomorphic to $\bold P^1_K-\beta^{-1}(0)$ and its special fiber is isomorphic to $\bold A^1_k$. 
Consider now the affine scheme $\Cal U=\text{Spec}\left(R[X,\beta(X)^{-1}]\right)$. We have a canonical injection $R[X]@>>>R[X,\beta(X)^{-1}]$ which induces an open immersion $\iota_0:\Cal U@>>>\Cal C_0$. Since $\beta(0)=0$, we have $X^{-1}\in R[X,\beta(X)^{-1}]$ and we can define an inclusion  $R[X_\infty,\gamma(X_\infty)^{-1}]@>>>R[X,\beta(X)^{-1}]$ by sending $X_\infty$ to $X^{-1}$. We then obtain a second open immersion $\iota_1:\Cal U@>>>\Cal D_1$. The morphisms $\iota_0$ and $\iota_1$ allow us to glue $\Cal C_1$ and $\Cal C_2$. We obtain a scheme $\Cal C=\Cal C_0\cup\Cal C_1\cong\bold P^1_{R}$, with $\Cal C_0\cap\Cal C_1=\Cal U$. Similarly, we can consider the affine scheme $\Cal V=\text{Spec}(R[T,T^{-1}])$ and the open immersions $\Cal V@>>>\Cal D_0$ and  $\Cal V@>>>\Cal D_1$, obtained, in terms of affine algebras, by taking the canonical injection $R[T]@>>>R[T,T^{-1}]$ and the inclusion $R[T_\infty]@>>>R[T,T^{-1}]$ which sends $T_\infty$ to $T^{-1}$. We then obtain a scheme $\Cal D=\Cal D_0\cup\Cal D_1\cong\bold P^1_{R}$, with $\Cal D_0\cap\Cal D_1=\Cal V$.  One easily checks that the morphisms $\Cal C_0@>>>\Cal D_0$ and $\Cal C_1@>>>\Cal D_1$ agree over $\Cal U$. In other words, we obtain a commutative diagram 
$$\CD
\Cal C_0@<\quad<<\Cal U@>>>\Cal C_1\\
@VVV @VVV @VVV\\
\Cal D_0@<\quad<<\Cal V@>>>\Cal D_1\endCD$$
which leads to a finite morphism 
$$\Cal C@>>>\Cal D$$
This is a semi-stable model for the cover over $R$ (since its generic fiber is isomorphic to $\beta$ over $K$) and we will call it the {\it fundamental model} associated to $\beta$. Lemma 3.1 implies that its specialization is purely inseparable. In particular, for a given $\lambda\in\Cal D(K)$, the elements of $\beta^{-1}(\lambda)$ all have the same specialization in the special fiber $\Cal C_k$ of $\Cal C$, so that this morphism will never separate the ramified fibers. A model $\Cal X@>>>\Cal Y$ for $\beta$ which separates the ramified fibers  will be obtained from   $\Cal C@>>>\Cal D$ after a finite number of blows-up. In order to get a minimal one, we may need to  blow-down the original irreducible component $\Cal C_k$ (and $\Cal D_k$ of the special fiber of $\Cal Y$). This last possibility will never occur. Indeed, this would mean that all the finite branch points have the same specialization in $D$, which is impossible, since we are assuming $\{0,1\}\in S$.

\specialhead \S4 Some complements on semi-stable models of genus zero\endspecialhead

 Before continuing, let's make some general remarks. First of all, suppose that $\Cal X/R$ is any  semi-stable model of a genus zero (pointed) curve  and denote by $\Cal X_k$ its special fiber. We classically define the {\it intersection graph} $\Gamma(\Cal X_k)$  as the abstract graph  whose vertices (resp.  edges) correspond to the irreducible components of $\Cal X_k$ (resp. the singular points of $\Cal X_k$). 

\midinsert
\epsfysize=2.8cm
\centerline{\epsfbox{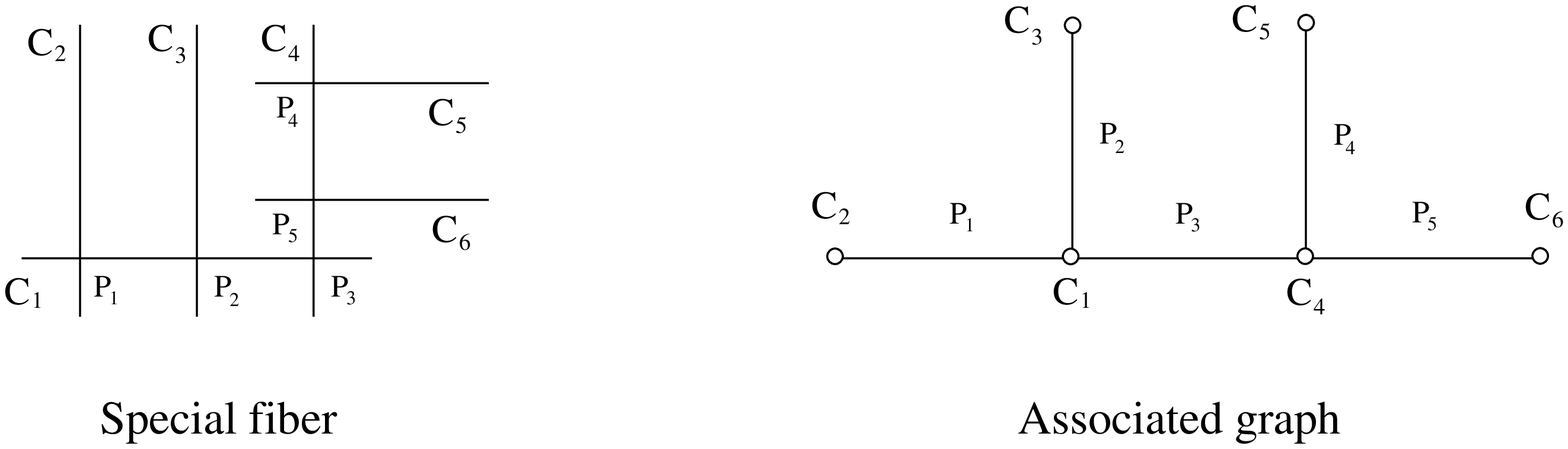}}
\endinsert

In this special case, $\Gamma(\Cal X_k)$ is a tree, endowed   with a {\it metric}, obtained by associating to any edge the thickness of the corresponding singular point. If $C_1$ and $C_2$ are two irreducible components of $\Cal X_k$, we can consider the {\it segment} $[C_1,C_2]$, which is the minimal (linear) subtree of $\Gamma(\Cal X_k)$ containing the vertices corresponding to $C_1$ and $C_2$. In particular, the distance $d(C_1,C_2)$ between $C_1$ and $C_2$ is defined as the sum of the lengths of the edges belonging to $[C_1,C_2]$. From a practical point of view, consider the model $\Cal X'$ of the projective line over $R$ obtained by blowing down all the irreducible components of $\Cal X_k$ different  from $C_1$ and $C_2$. Then, the special fiber of $\Cal X'$ is the union of $C_1$ and $C_2$, meeting at a unique singular point, of thickness $d(C_1,C_2)$.  
\midinsert
\epsfysize=2.6cm
\centerline{\epsfbox{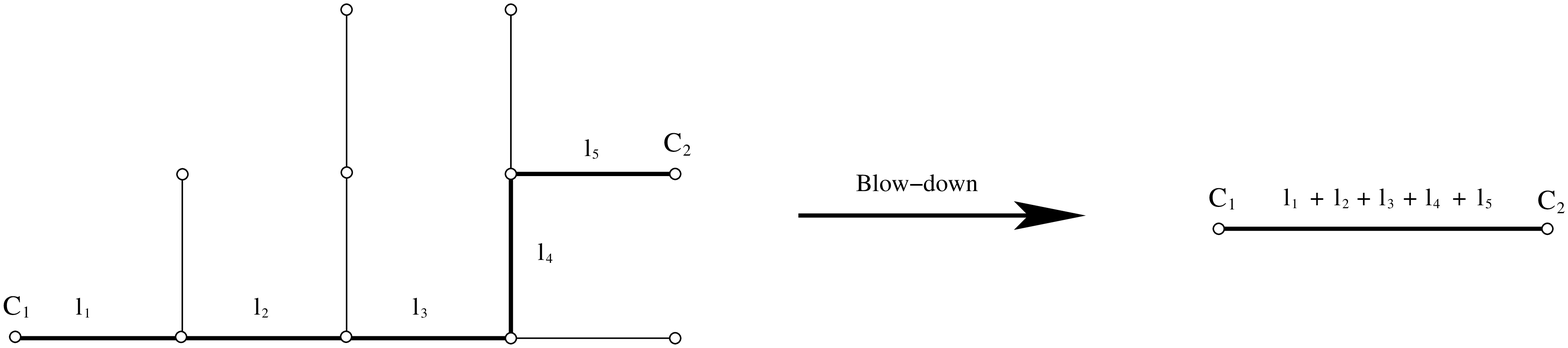}}
\endinsert

 Suppose now that $\Cal X@>>>\Cal Y$ is any semi-stable model over $R$ for a polynomial cover, obtained from its fundamental model $\Cal C@>>>\Cal D$ (cf. \S3) after a finite number of blows-up. Denote by $C$ (resp. by $D$) the irreducible component of $\Cal X_k$ (resp. of $\Cal Y_k$) corresponding to the special fiber of $\Cal C$ (resp. of $\Cal D$). Let $C_1\neq C$ be a tail of $\Cal C_k$ and denote by $D_1\neq D$  its image in $\Cal D_k$. If $C_\infty$ is the irreducible component of $\Cal X_k$ containing the specialization of the totally ramified point $\infty$, then we will suppose that the segments $[C,C_1]$ and $[C,C_\infty]$ have no common edges. Consider a point $P\in\Cal X(K)=\bold P^1(K)$ specializing in a non-singular point of $C_1$. Since $\Cal X(K)=\Cal C(K)$, the point $P$ defines a well defined element $x$ of $\Cal C(K)=\Cal C(R)$. Moreover, the above assumption on the relative positions of $C,C_0$ and $C_\infty$ imply that $x$ belongs to $\Cal C_0(R)=\bold A^1(R)$ (cf. the end of \S3) and thus, it can be assimilated to an element of $R$. Similarly the image of $P$ in $\Cal Y$ defines an element $\lambda$ of $\Cal D_0(R)=R$. Let $\pi$ and $\pi'$ be two elements of $\frak p$ such that $v(\pi)=d(C,C_1)$ and $v(\pi')=d(D,D_1)$. We then have a unique factorization
$$\beta(\pi X+x)-\lambda=\pi'\beta_0(X)\gamma_0(X)\tag{$*$}$$
where  $\beta(X)\in R[X]$ is the normalized model from which the cover $\Cal C@>>>\Cal D$ was constructed, $\beta_0(X)\in R[X]$ is monic and $\gamma_0(X)\in R[X]$ satisfies $\overline\beta(X)\in k^*$, i.e. $\gamma_0(X)\in R^*+\frak pR[X]$. Moreover, the finite morphism $\bold A^1_k@>>>\bold A^1_k$ obtained from $C_1@>>>D_1$ by removing the singular points is (isomorphic to the one) induced by the inclusion $k[T]@>>>k[X]$ which maps $T$ to $\overline\beta_0(X)$. In particular, $C_1$ will be the only irreducible component of $\Cal X_k$ lying above $D_1$ if and only if $\gamma_0(X)$ has degree $0$, i.e. if $\gamma_0(X)=1$. In this case, we get the relation $d(D,D_1)=pd(C,C_1)$.

\specialhead \S5 Simple reduction of the branch locus\endspecialhead

 We will now put some conditions on the reduction type of the branch locus $B=S\cup\{\infty\}$ of a polynomial  cover over $K$. The curve $\Cal D$ of the fundamental model associated to it (cf. \S2) is a model of the projective line over $R$ such that no point of $S$ specialize to $\infty$ and at least two points have distinct specializations. These two conditions uniquely determine $\Cal D$, up to $R$-isomorphism. The stable model $\Cal B$ (over $R$) for the pointed curve $(\bold P^1_K,B)$ is obtained from $\Cal D$ after a finite number of blows-up.

\definition{5.1 Definition} Consider the stable (separating) model $\Cal B$ associated to a pointed curve $(\bold P^1_K,B)$, where $B=S\cup\{\infty\}$. Denote by $D$ the irreducible component of its special fiber containing the specialization of the point $\infty$. An element $\lambda\in S$ is {\it ordinary} if it specializes in $D$. If $\Cal B$ has bad reduction,  an irreducible component  of $\Cal B_k$  is {\it simple}  if it only meets $D$ and if there are exactly two elements of $S$ specializing in it. We will say that  $B$  has {\it simple reduction} if  any tail of $\Cal B_k$ (if there are any) is simple (see the following picture).

\vskip.6cm
 
\epsfysize=2.4cm

\centerline{\epsfbox{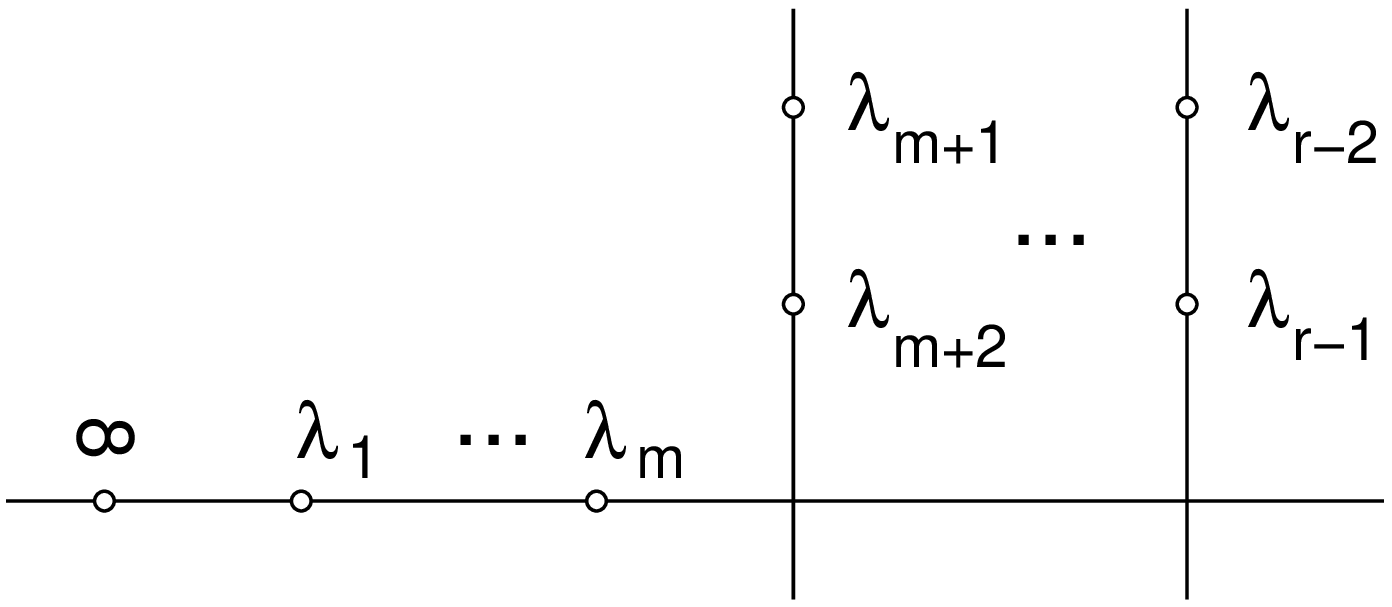}}

\vskip.6cm 

\enddefinition

Remark that if the cardinality of $B$ is less than or equal to four, then the pointed curve $(\bold P^1_K,B)$ has automatically simple reduction.

\specialhead \S6 Semi-stable model separating the fiber above an ordinary branch point\endspecialhead

 We will now describe the stable model for a polynomial cover of degree $p$ dominating its fundamental model $\Cal C@>>>\Cal D$ and separating the ramified fiber above an ordinary branch point. 

\proclaim{6.1 Theorem} Let $\beta:\bold P^1_K@>>>\bold P^1_K$ be a polynomial cover of degree $p$. Denote by $B=S\cup\{\infty\}$ its branch locus and suppose that $\lambda\in S$ is ordinary. Then, the minimal semi-stable model  $\Cal X@>>>\Cal Y$ of the cover  which dominates the fundamental model $\Cal C@>>>\Cal D$ and separates the fiber above $\lambda$ has the following description:

\roster
\item The special fiber of $\Cal X$  is the union of two projective lines $C$ and $C'$ meeting at a unique singular point of thickness $(n_\lambda-1)^{-1}$ where $n_\lambda$ is the cardinality of the fiber $\beta^{-1}(\lambda)$.
\item The special fiber of $\Cal Y$ is the union of two projective lines  $D$ and $D'$ meeting at a unique singular point of thickness $p(n_\lambda-1)^{-1}$.
\item  The morphism $C@>>>D$ is purely inseparable, while $C'@>>>D'$ is generically \'etale, unramified outside two points, wildly ramified above one of them and tamely ramified above the other.
\endroster

\vskip.6cm
 
\epsfysize=2cm

\centerline{\epsfbox{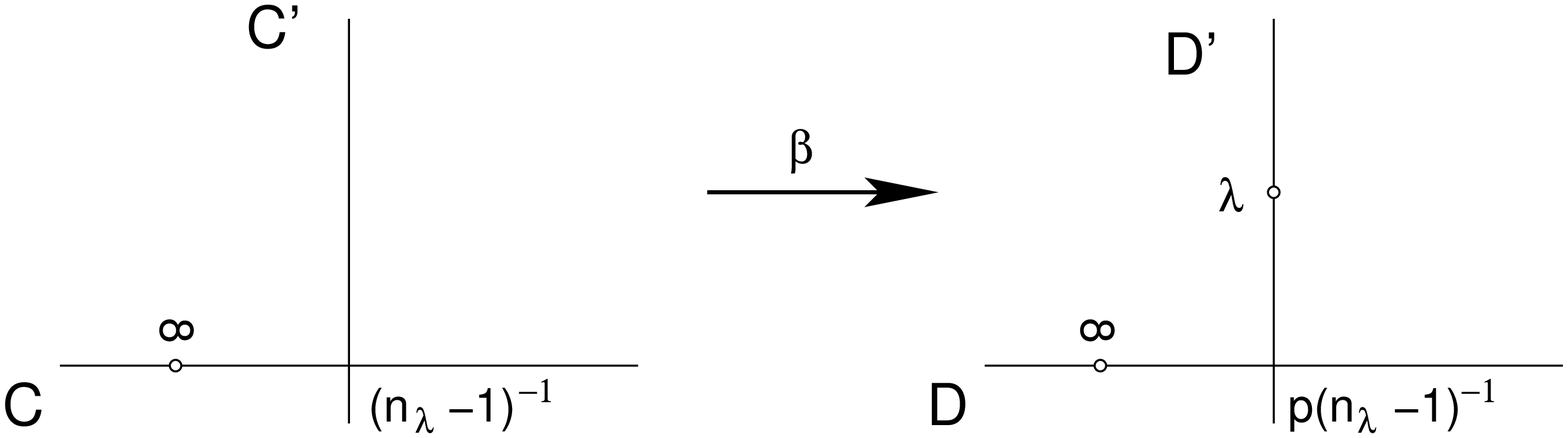}}

\vskip.6cm 

\endproclaim

\demo{Proof} Without any loss of generality, we can assume that $\lambda=0$  and $\beta(0)=0$. Let $\nu$ be the rational number defined by 
$$\nu=\text{Min}\{v(x)\,\,|\,\,x\in\beta^{-1}(0)\}$$ 
and consider the curve $\Cal C'$ obtained from $\Cal C$ after a blow-up at the origin, of thickness $\nu$. Let $\pi$ be an element of $\frak p$ such that $v(\pi)=\nu$. By construction, we have the expression
$$\beta(\pi X)=\pi^p\prod_{x\in\beta^{-1}(0)}(X-\pi^{-1}x)^{e_x}\in R[X]$$
with $\pi^{-1}x\in R$ for any $x\in\beta^{-1}(0)$. Following the notation introduced at the end of \S4, we obtain  $\beta_0(X)=\pi^{-p}\beta(\pi X)$ and $\gamma_0(X)=1$. In particular, we have a model $\Cal C'@>>>\Cal D'$ of the cover $\beta$, where $\Cal D'$ is the curve obtained from $\Cal D$ after a blow-up at the origin, of thickness $p\nu$. We just have to prove that it separates the fiber above $0$. From a practical point of view, we must show that two distinct roots of the polynomial $\beta_0(X)$ have  distinct stpecializations. First of all, if $\lambda\in S-\{0\}$ then $\beta^{-1}(\lambda)\subset R^*$. Indeed, we already have $\beta^{-1}(\lambda)\subset R$. If there were $x\in\beta^{-1}(\lambda)$ with positive valuation, we would obtain $\overline\lambda=\overline\beta(x)=\overline x^p=0$, which is impossible since we are assuming that $0\in S$ is ordinary. We have the following expression for the derivative of $\beta(X)$:
$$\beta'(X)=p\prod_{x\in\beta^{-1}(S)}(X-x)^{e_x-1}$$
We then obtain $\beta_0'(X)=\pi^{1-p}\beta'(\pi X)$, which leads to
$$\beta_0'(X)=p\pi^{1-n_0}\prod_{x\in\beta^{-1}(0)}(X-\pi^{-1}x)^{e_x-1}\prod_{x\in\beta^{-1}(S-\{0\})}(\pi X-x)^{e_x-1}$$
where $n_0$ is the cardinality of the fiber $\beta^{-1}(0)$. We have $\beta'_0(X)\in R[X]$, which directly implies the inequality $(n_0-1)\nu\leq1$. A strict inequality would give $\overline\beta_0'(X)=0$ and thus $\overline\beta_0(X)=X^p$ (recall that $\beta_0(X)$ is monic, of degree $p$ and satisfies $\beta_0(0)=0$). This is impossible since, by construction,  there exists a root of $\beta_0(X)$ belonging to $R^*$. We then have $\nu={1\over n_0-1}$. Suppose now that the elements $x_1,\dots,x_s\in\beta_0^{-1}(0)$ specialize to the same point $t$ of $C'$ and denote by $e_1,\dots,e_s$ their ramification indices. We have $s<n_\lambda$, since there are at least two distinct roots of $\overline\beta_0(X)$. The monomial $X-t$ will appear  with a multiplicity $e=e_1+\dots+e_s<p$ in the factorization of $\overline\beta_0(X)$. In particular, the element $t$ will be a root of $\overline\beta_0'(X)$ of multiplicity $e-1$. But the above expression of $\beta_0'(X)$ implies that this multiplicity is equal to $e-s$, and so we have $s=1$, i.e. the model $\Cal C'@>>>\Cal D'$ separates the fiber above $0$.\qed
\enddemo

\proclaim{6.2 Corollary} The notation and hypothesis being as in theorem 1, there exists a smooth model of the cover (over $R$) which separates the fibers above $\infty$ and $\lambda$. In other words, the pointed curve $(\bold P^1_K,\beta^{-1}(\{\infty,\lambda\}))$ has good reduction.
\endproclaim

\demo{Proof} The desired smooth model is obtained  from $\Cal C'@>>>\Cal D'$ by  blowing down the irreducible components $C$ and $D$.\qed\enddemo

\proclaim{6.3 Corollary} With the above assumptions, suppose that the branch locus $B=\{\infty,\lambda_1,\dots,\lambda_{r-1}\}$ of the cover $\beta$ has good reduction and that $r\geq3$. Then, the minimal stable model $\Cal X@>>>\Cal Y$ separating the ramified fibers has the following description: 

\roster 
\item $\Cal X_k$ (resp. $\Cal Y_k$) is the union of $r$ projective lines $C,C_1,\dots,C_{r-1}$ (resp. $D,D_1,\dots,D_{r-1}$). 
\item  For any $i\in\{1,\dots,r-1\}$, the point $\lambda_i$ specializes in $D_i$.
\item  The irreducible component $C_i$ (resp. $D_i$) only meets $C$ (resp. $D$). The corresponding  singular point has thickness $(n_i-1)^{-1}$ (resp. $p(n_i-1)^{-1}$), where $n_i$ denotes the cardinality of the fiber above $\lambda_i$.

\item The morphism $C@>>>D$ is purely inseparable  and, for any $i\in\{1,\dots,r-1\}$, the cover $C_i@>>>D_i$ is generically \'etale, unramified outside two points, wildly ramified above one of them (the singular point) and tamely ramified above the other (the specialization of $\lambda_i$).
\endroster
\vskip.6cm
 
\epsfysize=2.4cm

\centerline{\epsfbox{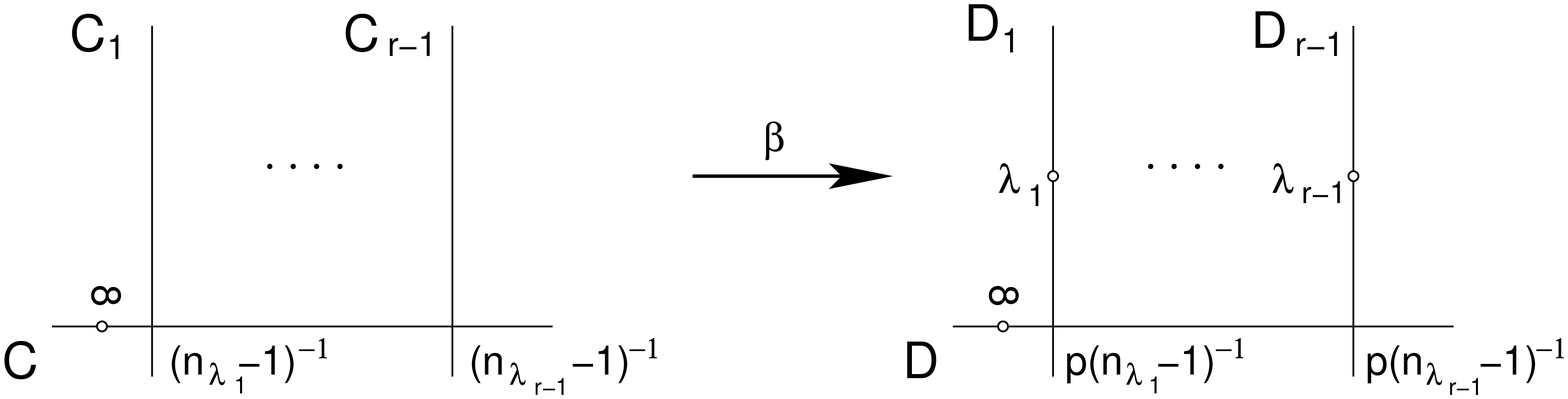}}

\vskip.6cm 
\endproclaim

\demo{Proof} It suffice to repeat the construction in the proof of Theorem 6.1 for all the elements of $S$.
\qed\enddemo

\specialhead \S7 Semi-stable model separating the fibers above a simple tail\endspecialhead

 We will now study the case of simple reduction of the branch locus $B=S\cup\{\infty\}$ of the cover $\beta$. Let $r$ be the cardinality of $B$. We can assume $r>3$, since for $r\leq 3$ the pointed curve $(\bold P^1,B)$ will always have good reduction. If $\lambda\in S$ is an ordinary branch point, then the construction in the proof of theorem 1 leads to a stable model  which separates the fiber above $\lambda$. Suppose that $\lambda,\lambda'\in S$ belong to a simple tail of the minimal stable model $\Cal B$ of $B$. In other words, viewing $S$ as a subset of $R=\Cal C_0(R)$, we have $v(\lambda-\lambda')>0$ and $v(\lambda-\lambda'')=v(\lambda'-\lambda'')=0$ for any $\lambda''\in S-\{\lambda,\lambda'\}$. The positive rational number $\epsilon=v(\lambda-\lambda')$ is the thickness of the singular point connecting the simple tail with the rest of $\Cal B_k$. Since the reduction of the model $\Cal C@>>>\Cal D$ is purely inseparable, all the elements of $\beta^{-1}(\{\lambda,\lambda'\})$ will have the same specialization. 

\subsubhead a) The first separating blow-up\endsubsubhead We will start the construction of the stable model separating the fibers above $\lambda$ and $\lambda'$ by considering the minimal model $\Cal C'$ dominating $\Cal C$ such that  all the elements of  $\beta^{-1}(\{\lambda,\lambda'\})$ specialize in the (smooth locus of the) same irreducible component and at least two of them have distinct specializations. 
In order to obtain an more explicit description, we can first of all reduce to the case $\lambda'=0$ and $\beta(0)=0$. Consider the rational number $\nu$ defined by 
$$\nu=\text{Min}\left\{v(x-y)\,\,|\,\,x\in\beta^{-1}(\{0,\lambda\})\right\}$$
We have $\nu>0$, since $\overline\lambda=0$ and $\overline\beta(X)=X^p$. Moreover, $v(\lambda)\geq p\nu$. The curve $\Cal C'$ is obtained from $\Cal C$  after a blow-up at the origin, of thickness $\nu$. Let $\pi$ be an element of $\frak p$ such that $v(\pi)=\nu$. The relation $(*)$ at the end of \S4 becomes $\beta(\pi X)=\pi^p\beta_0(X)\gamma_0(X)$, with  $\gamma_0(X)=1$ and
$$\beta_0(X)=\pi^{-p}\beta(\pi X)=\prod_{x\in\beta^{-1}(0)}(X-\pi^{-1}x)^{e_x}=\pi^{-p}\lambda+\prod_{x\in\beta^{-1}(\lambda)}(X-\pi^{-1}x)^{e_x}$$
We then get a model $\Cal C'@>>>\Cal D'$ for $\beta$, where $\Cal D'$ is the curve obtained from $\Cal D$ after a blow-up at the origin, of thickness $p\nu$ (see the following picture).

\vskip.6cm
 
\epsfysize=2.2cm

\centerline{\epsfbox{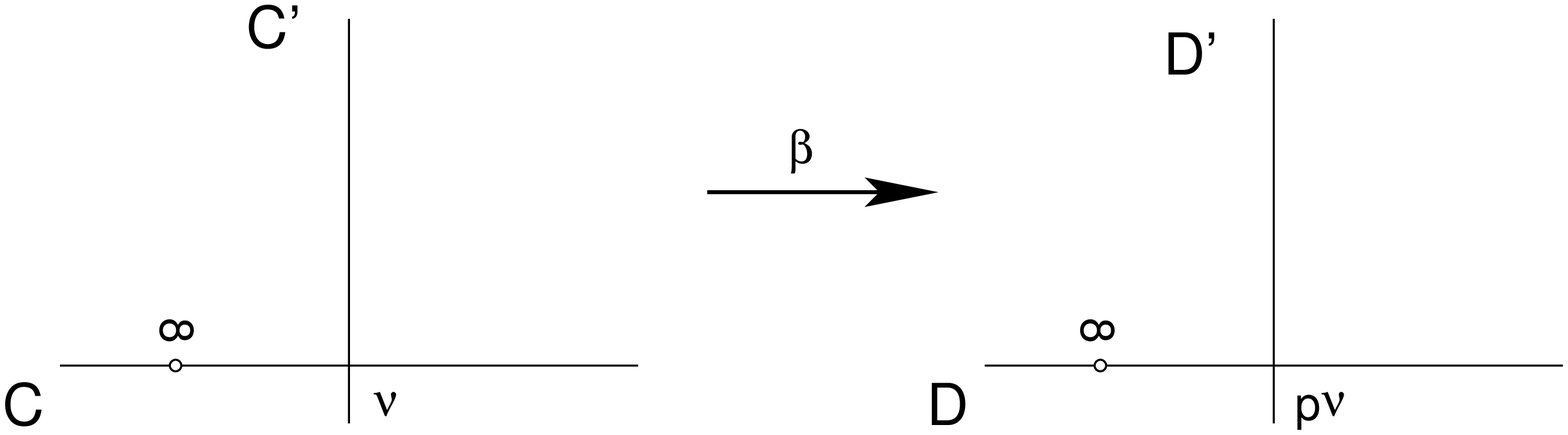}}

\vskip.6cm 

\noindent The derivative of $\beta_0(X)$ can be expressed as 
$$\beta_0'(X)=p\pi^{p+1-n_0-n_\lambda}\prod_{x\in\beta^{-1}(\{0,\lambda\})}(X-\pi^{-1}x)^{e_x-1}\prod_{x\in\beta^{-1}(S-\{0,\lambda\})}(\pi X-x)^{e_x-1}$$
where, for any $\lambda\in S$,  $n_\lambda$ is the cardinality of the fiber $\beta^{-1}(\lambda)$. In particular, since $\beta_0'(X)\in R[X]$ and $\pi^{-1}x\in R$ for any $x\in\beta^{-1}(\{0,\lambda\})$, we obtain  the relation 
$$\nu\leq{1\over n_0+n_\lambda-p-1}\tag{$**$}$$ 

\subsubhead b) First case: ``far'' points\endsubsubhead  Let's start by supposing that this last inequality is strict. This directly implies that $\overline\beta_0(X)=X^p$. Now, by construction, the polynomial $h(X)=\beta_0(X)(\beta_0(X)-\pi^{-p}\lambda)$ has at least two roots having distinct specializations, so that we obtain   
$$v(\lambda)=p\nu$$
 Indeed, the relation $v(\lambda)>p\nu$  would give $\overline h(X)=X^{2p}$, which has only one root. In particular, this frist blow-up separates the branch points $0$ and $\lambda$ but does not separate the corresponding fibers.

\vskip.6cm
 
\epsfysize=2.4cm

\centerline{\epsfbox{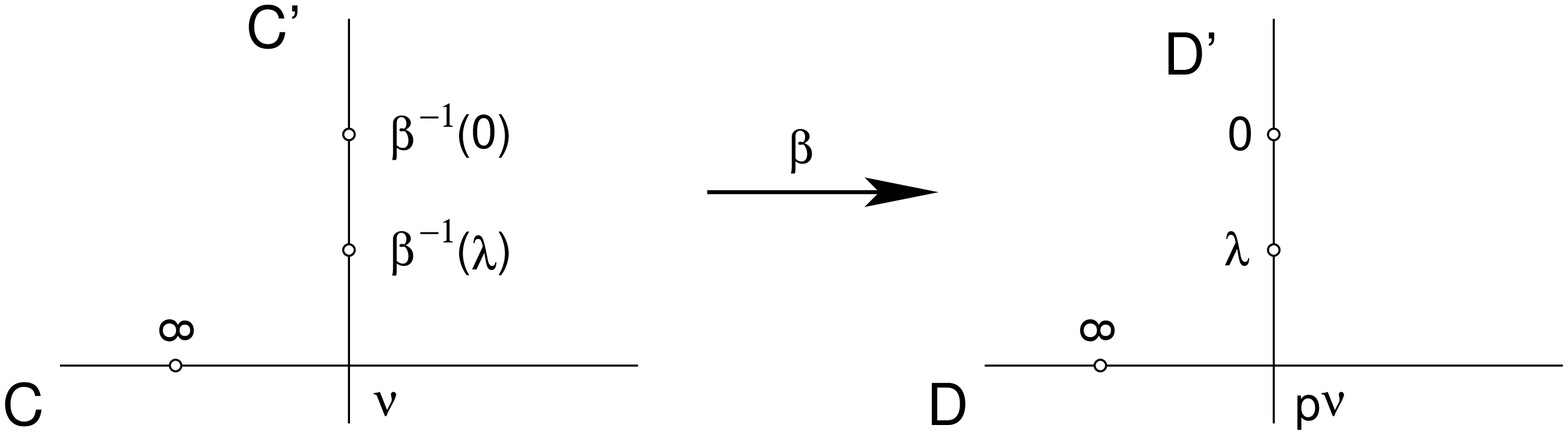}}

\vskip.6cm 

\noindent  The above inequality now reads as  
$$v(\lambda)<{p\over n_0+n_\lambda-p-1}$$
In order to obtain a model separating the fiber above $0$, let $\nu_0$ be the positive integer defined by 
$$\nu_0=\text{Min}\left\{v(x)\,\,|\,\,x\in\beta_0^{-1}(0)\right\}$$
We clearly have $\nu+\nu_0=\text{Min}\left\{v(x)\,\,|\,\,x\in\beta^{-1}(0)\right\}$. Consider the curve $\Cal C''$ obtained from $\Cal C'$ after a blow-up at the origin of the irreducible component $C'$, of thickness $\nu_0$. If $C''$ denotes the new irreducible component of $\Cal C''_k$, then we have $d(C'',C)=\nu+\nu_0$. Let $\pi_0\in\frak p$ such that $v(\pi_0)=\nu_0$ and set $\pi_1=\pi\pi_0$ (as in the previous paragraph, $\pi$ has valuation $\nu$). We then obtain the expression $\beta(\pi_1X)=\pi^p\beta_0(\pi_0X)=\pi_1^p\beta_1(X)\gamma_1(X)$, with  $\gamma_1(X)=1$ and 
$$\aligned\beta_1(X)&=\pi_1^{-p}\beta(\pi_1X)=\pi^{-p}\beta_0(\pi_0 X)=\\
&=\prod_{x\in\beta^{-1}(0)}(X-\pi_1^{-1}x)^{e_x}=\prod_{x\in\beta_0^{-1}(0)}(X-\pi_0^{-1}x)^{e_x}\endaligned$$
By construction, at least two roots of $\beta_1(X)$ have distinct specializations. If $\Cal D''$ is the curve obtained from $\Cal D'$ after a blow-up at the origin of $D'$, of thickness $p\nu_0$, we then have a model $\Cal C''@>>>\Cal D''$ of the cover $\beta$. Setting $u=n_0+n_\lambda-p-1$, $u_0=n_0-1$ and $\delta=1-u\nu-u_0\nu_0$,
we have the expression 
$$\aligned\beta_1'(X)=&p\pi^u\pi_0^{u_0}\prod_{x\in\beta^{-1}(0)}(X-\pi_1^{-1}x)^{e_x-1}\prod_{x\in\beta^{-1}(\lambda)}(\pi_0X-\pi^{-1}x)^{e_x-1}\times\\
&\times\prod_{x\in\beta^{-1}(S-\{0,\lambda\})}(\pi_1X-x)^{e_x-1}\endaligned$$
which implies that $\delta=0$ (otherwise, we would obtain $\overline\beta_1(X)=X^p$ and all the roots of $\beta_1(X)$ would have the same specialization). Since $v(\lambda)=p\nu$, we obtain
$$\nu_0={p-(n_0+n_\lambda-p-1)v(\lambda)\over p(n_0-1)}<\frac1{n_0-1}$$
Proceeding exactly as in the end of the proof of theorem 1, we finally check that this model separates the fiber above $0$, i.e. that any two distinct roots of $\beta_1(X)$ have distinct specializations. The same procedure leads to a stable model separating the fiber above $\lambda$. Summarizing, we have just proved the following result:

\proclaim{7.1 Theorem} Suppose that the branch locus $B=S\cup\{\infty\}$ of the polynomial cover $\beta$ has bad reduction and that $D'$ is a simple tail of the special fiber of the stable model $\Cal B$ associated to the pointed curve $(\bold P^1,B)$. Denote by  $\lambda_1,\lambda_2\in S$  the two branch points specializing in $D'$ and let $\epsilon=v(\lambda_1-\lambda_2)$ be the thickness of the corresponding singular point of $\Cal B_k$. For any $\lambda\in S$, denote by $n_\lambda$ the cardinality of the fiber $\beta^{-1}(\lambda)$. If the inequality 
$$\epsilon<{p\over n_{\lambda_1}+n_{\lambda_2}-p-1}$$
 holds, then, the minimal semi-stable model $\Cal X@>>>\Cal Y$ of $\beta$ dominating $\Cal C@>>>\Cal D$ and   separating the fibers above  $\lambda_1$ and $\lambda_2$ has the following description: 

\roster

\item The curve $\Cal X_k$ (resp. $\Cal Y_k$) is the union of four projective lines $C,C',C_1$ and $C_2$ (resp. $D,D',D_1$ and $D_2$).

\item The specializations of the points $\infty,\lambda_1$ and $\lambda_2$ belong respectively to $D,D_1$ and $D_2$.

\item The irreducible component $C$ (resp. $D$) only meets $C'$ (resp. $D'$), at a singular point of thickness $\frac\epsilon p$ (resp. of thickness $\epsilon$).

\item For any $i\in\{1,2\}$, the irreducible component $C_i$ (resp. $D_i$) only meets $C'$ (resp. $D'$), at a singular point of thickness $\nu_i={p-(n_{\lambda_1}+n_{\lambda_2}-p-1)\epsilon\over p(n_{\lambda_i}-1)}$ (resp. of thickness $p\nu_i={p-(n_{\lambda_1}+n_{\lambda_2}-p-1)\epsilon\over n_{\lambda_i}-1}$).

\item The morphisms $C@>>>D$ and $C'@>>>D'$ are purely inseparable, while the covers  $C_1@>>>D_1$ and  $C_2@>>>D_2$ are generically \'etale, unramified outside two points, wildly ramified above one of them  and tamely ramified above the other.

\endroster
\vskip.6cm
 
\epsfysize=2.4cm

\centerline{\epsfbox{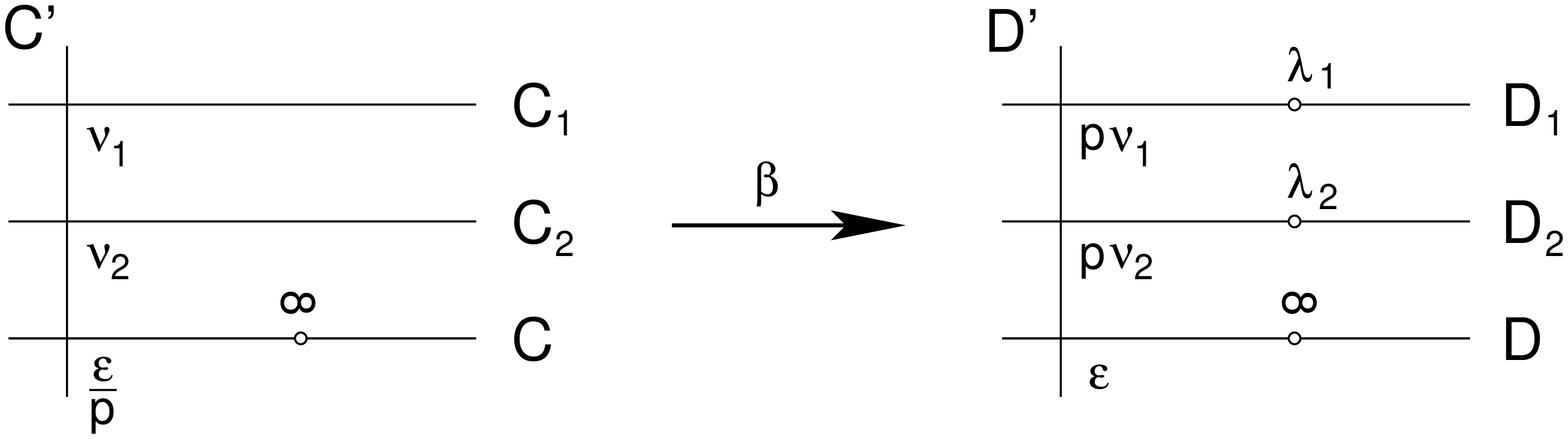}}

\vskip.6cm 
\endproclaim 

\proclaim{7.2 Corollary} The notation and hypothesis being as above, if $\epsilon<{p\over n_{\lambda_1}+ n_{\lambda_2}-p-1}$ then, for any  $i\in\{1,2\}$, the pointed curve $(\bold P^1,\beta^{-1}(\{\infty,\lambda_i\}))$ has good reduction.\endproclaim

\demo{Proof} It suffices to take the smooth $R$-curve obtained from $\Cal X$ by blowing down the irreducible components $C$, $C'$ and $C_{2-i}$ of $\Cal C_k$.\qed\enddemo

\subsubhead c) Second case: the critical distance\endsubsubhead  Keeping the above notation and hypothesis, suppose now that the inequality $(**)$ at the end  of \S7a is in fact an equality, i.e. that  $\nu={1\over n_0+n_\lambda-p-1}$, so that $\overline\beta_0'(X)\neq 0$. In particular, the morphism $C'@>>>D'$ is generically \'etale. We already know that $v(\lambda)\geq p\nu$. In this paragraph, we will assume that this last inequality is an equality. As before, the cover $\Cal C'@>>>\Cal D'$ will separate the branch points $0$ and $\lambda$. Moreover, proceeding as at the end of the proof of Theorem 6.1, we can easily  show that this model also separates the fibers above these two points. In particular, we have the following result:

\proclaim{7.3 Theorem} The notation and hypothesis being as in Theorem 7.1, suppose that
$$\epsilon={p\over n_{\lambda_1}+n_{\lambda_2}-p-1}$$
 Then, the minimal semi-stable model $\Cal X@>>>\Cal Y$ of $\beta$ dominating $\Cal C@>>>\Cal D$ and separating  the fibers above $\lambda_1$ and $\lambda_2$ has the following description: 

\roster

\item The curve $\Cal X_k$ (resp. $\Cal Y_k$) is the union of two projective lines $C$ and $C'$ (resp. $D$ and $D'$) meeting at a singular point of thickness $\frac\epsilon p$ (resp. of thickness $\epsilon$).

\item The specializations of the points $\lambda_1$ and $\lambda_2$ belong $D'$ and $\infty$ specializes in $D$.

\item The morphism $C@>>>D$ is purely inseparable while the cover $C'@>>>D'$ is generically \'etale, unramified outside three points, wildly ramified above one of them and tamely ramified above the others.

\endroster
\vskip.6cm
 
\epsfysize=2.4cm

\centerline{\epsfbox{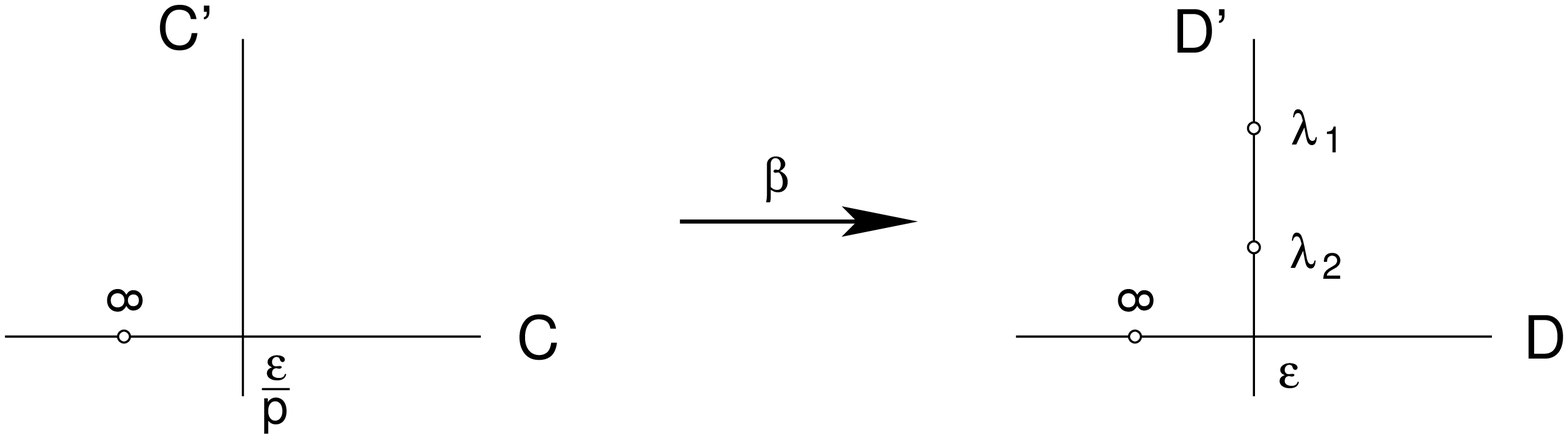}}

\vskip.6cm 
\endproclaim 

\proclaim{7.4 Corollary} The notation and hypothesis being as above, if $\epsilon={p\over n_{\lambda_1}+n_{\lambda_2}-p-1}$ then the pointed curve $(\bold P^1,\beta^{-1}(\{\infty,\lambda_1,\lambda_2\}))$ has good reduction.\endproclaim

\demo{Proof} It suffices to take the smooth $R$-curve obtained from $\Cal X$ by blowing down the irreducible component $C$ of its special fiber.\qed\enddemo
 
\subsubhead d) Third case: ``near'' points\endsubsubhead  We now have to study the most critical case, that is, when $\nu={1\over n_0+n_\lambda-p-1}$ and $v(\lambda)>p\nu$. This situation occurs if and only if the model $\Cal C'@>>>\Cal D'$ introduced at the beginning of this section does not separate the elements $0$ and $\lambda$ of the branch locus $B$ of the cover. As in the previous paragraph, the derivative of $\overline\beta_0(X)$ will not identically vanish, so that the cover $C'@>>>D'$, which is (isomorphic to the one) induced by the polynomial $\overline\beta_0(X)\in k[X]$,  is generically \'etale, unramified outside two points, wildly ramified above one of them (the intersection with $D$) and tamely ramified above the other . There exist at least two roots of $\beta_0(X)$ having distinct specializations. Indeed, the contrary would give $\overline\beta_0(X)=X^p$, so that the cover $C'@>>>D'$ would be purely inseparable, which is a contraddiction. The points $0$ and $\lambda$ are the only elements of $B$ specializing in $C'$ and their fibers with respect to $\beta$ can be assimilated to the fibers $\beta_0^{-1}(0)$ and $\beta_0^{-1}(\lambda_0)$, where $\lambda_0=\pi^{-p}\lambda\in\frak p$. Set 
$$\overline\beta_0(X)=\prod_{i=1}^s(X-w_i)^{d_i}$$
with $s>1$, $d_1+\dots+d_s=p$ and $w_i\neq w_j$ for any $i\neq j$. We then obtain two partitions 
$\beta_0^{-1}(0)=S_{1,0}\cup\dots\cup S_{s,0}$ 
and $\beta_0^{-1}(\lambda_0)=S_{1,\lambda_0}\cup\dots\cup S_{s,\lambda_0}$,
 where we have set, for a fixed $t\in\frak p$,
$$S_{i,t}=\{x\in\beta_0^{-1}(t)\,\,|\,\,\overline x=w_i\}$$
For any $i\in\{1,\dots,s\}$ we then have the identity 
$$d_i=\sum_{x\in S_{i,t}}e_x$$
 where $e_x$ is the multiplicity of the root $x$ of the polynomial $\beta_{0}(X)-t$. Since the cover $C'@>>>D'$ is generically \'etale and unramified outside the set $\{0,\infty\}$, we have 
$$\overline\beta_0'(X)=c\prod_{i=1}^s(X-w_i)^{d_i-1}$$
with $c\in k^*$. On the other hand, the expression of the derivative of $\beta_0(X)$ given previously leads to  
$$\overline\beta_0'(X)=c\prod_{i=1}^s(X-w_i)^{2d_i-n_{i,0}-n_{i,\lambda}}$$
where $n_{i,0}$ (resp. $n_{i,\lambda}$) denotes the cardinality of $S_{i,0}$ (resp. of $S_{i,\lambda_0}$). Combining these two expressions we finally obtain the identity $$n_{i,0}+n_{i,\lambda}=d_i+1$$
which holds for any $i\in\{1,\dots,s\}$. In order to completely separate the ramified fibers, we need to blow-up some projective lines at the points $w_1,\dots,w_s$ belonging to the irreducible component $C'$ of $\Cal C'$. In other words, we have to separate the elements of $S_{i,0}\cup S_{i,\lambda_0}$, for any $i\in\{1,\dots,s\}$. Let's start with $i=1$: since $0\in S$ is a branch point and $\beta_0(0)=\beta(0)=0$, we can assume, without any loss of generality,  $w_1=0$. Set 
$$\nu_1=\text{Min}\left\{v(x-y)\,\,|\,\,x,y\in S_{1,0}\cup S_{1,\lambda_0},\,\, x\neq y\right\}>0$$
 and consider the $R$-curve $\Cal C''$, obtained from $\Cal C'$ by  blowing-up a projective line $C_1$  at $w_1$, of thickness $\nu_1$. Let $\pi_2\in\frak p$  such that $v(\pi_2)=\nu_1$ and set $\pi_3=\pi\pi_2$. We then have the identity
$$\beta(\pi_3X)=\pi^p\beta_0(\pi_2X)=\pi_3^p\beta_1(X)\gamma_1(X)$$
where
$$\beta_1(X)=\prod_{x\in S_{1,0}}(X-\pi_2^{-1}x)^{e_x}$$
and
$$\gamma_1(X)=\prod_{x\in\beta_0^{-1}(0)-S_{1,0}}(\pi_2X-x)^{e_x}$$
Since we are assuming $w_1=0$, an element $x\in\beta_0^{-1}(0)$ belongs to $S_{1,0}$ if and only if $v(x)>0$. In particular, $\overline\gamma_1(X)=c\in k^*$. We then obtain a $R$-morphism $\Cal C''@>>>\Cal D''$ (which is not finite), where $\Cal D''$ is the $R$-curve obtained from $\Cal D'$ by blowing-up a projective line $D_1$ at the origin of $D'$, of thickness $\epsilon_1=p\nu_1$. For any  $x\in S_{i,\lambda_0}$, we have $v(x)\geq\nu_1$, from which we easily deduce the inequality $v(\lambda_0)\geq\epsilon_1$, i.e.  $v(\lambda)\geq\epsilon_0+\epsilon_1$, with $\epsilon_0=p\nu$.
The specialized morphism $C_1@>>>D_1$ has degree $d_1$ and is (isomorphic to the cover $\bold P^1_k@>>>\bold P^1_k$) induced by the polynomial $\overline\beta_1(X)$. The integer $d_1$ being strictly less than $p$, the derivative $\overline\beta_1'(X)$ does not identically vanish. More explicitely, using once again the expression of the derivative of $\beta_0(X)$, we obtain
$$\overline\beta_1'(X)=u\prod_{x\in S_{i,0}\cup S_{i,\lambda_0}}\left(X-\overline{\pi_2^{-1}x}\right)^{e_x-1}$$
with $u\in k^*$. The strict inequality $v(\lambda)>\epsilon_0+\epsilon_1$ would imply that the cover $\bold P^1_k@>>>\bold P^1_k$ induced by $\overline\beta_1(X)$ is tame, of degree $d_1$ and unramified outside $0$ and $\infty$, and thus $\overline\beta_1(X)=uX^{d_1}$. In particular, all the elements of  $S_{1,0}$ and $S_{1,\lambda_0}$, would specialize to the same element of $C_1$. This is impossible,  since $\Cal C''$ is the minimal model dominating $\Cal C'$ and separating at least two elements of  $S_{1,0}\cup S_{1,\lambda_0}$.  We then have the identity 
$$v(\lambda)=\epsilon_0+\epsilon_1={p\over n_0+n_\lambda-p-1}+d_1\nu_1$$
 which implies  that the model $\Cal D''$ separates the branch points $0$ and $\lambda$. Then, using the above expression for the derivative of $\overline\beta_1(X)$ and the same arguments at the end of the proof of Theorem 6.1, one easily shows that $\Cal C''$ separates the whole $S_{1,0}\cup S_{1,\lambda_0}$. Iterating this construction for all the roots of $\overline\beta_0(X)$, we obtain a semi-stable model $\Cal X@>>>\Cal Y$ (this time the morphism is finite) which dominates $\Cal C@>>>\Cal D$ and separates the fiber above $0$ and $\lambda$. Summarizing, we have proved the following result:

\proclaim{7.5 Theorem} The notation and hypothesis being as in theorem 2, suppose that
$$\epsilon>{p\over n_{\lambda_1}+n_{\lambda_2}-p-1}$$
Then, the special fiber of the minimal stable model $\Cal X@>>>\Cal Y$ of the cover $\beta$ dominating $\Cal C@>>>\Cal D$ and separating  the fibers above $\lambda_1$ and $\lambda_2$  has the following description: 

\roster

\item The curve $\Cal Y_k$  is the union of three projective lines $D$, $D'$ and $D_1$.

\item The curve $D$ (resp. $D_1$) meets $D'$ at a singular point of thickness $\epsilon_0={p\over n_{\lambda_1}+n_{\lambda_2}-p-1}$ (resp. of thickness $\epsilon_1=v(\lambda)-{p\over n_{\lambda_1}+n_{\lambda_2}-p-1}$).

\item The specializations of the points $\lambda_1$ and $\lambda_2$ belong $D_1$ and $\infty$ specializes in $D$.

\item The curve $\Cal X_k$  has $s+2\geq 4$ irreducible  components $C,C',C_1,\dots,C_s$.

\item The curve $C$ meets $C'$ at a singular point of thickness $\frac{\epsilon_0}p$ .

\item  There exist two partitions 
$$\beta^{-1}(\lambda_1)=S_{1,1}\cup\dots\cup S_{s,1}\qquad\text{and}\qquad \beta^{-1}(\lambda_2)=S_{1,2}\cup\dots\cup S_{s,2}$$
 satisfying the identities
$$\sum_{x\in S_{i,1}}e_x=\sum_{x\in S_{i,2}}e_x=d_i\qquad\text{and}\qquad n_{i,1}+n_{i,2}=d_i+1$$
 where $e_x$ denotes the ramification index of $\beta$ at a point $x\in\bold P^1(K)$and, for any $i\in\{1,\dots,s\}$ and 
 any $j\in\{1,2\}$,  $n_{i,j}$ is the cardnality of the set $S_{i,j}$.

\item For any $i\in\{1,\dots,s\}$, the curve $C_i$ meets $C$ at a unique singular point, of thickness $\frac{\epsilon_1}{d_i}$ and the elements of $S_{i,1}\cup S_{i,2}$ specialize to pairwise distinct points of $C_i$.

\item The morphism $C@>>>D$ is purely inseparable and the cover  $C'@>>>D'$ is generically \'etale, unramified  outside two points, wildly ramified above one of them (the intersection with $D$) and tamely ramified above the other (the intersection with $D_1$). For any $i\in\{1,\dots,s\}$, the cover $C_i@>>>D_1$ is generically \'etale, of degree $d_i$, unramified outside three points (the specializations of $\lambda_1$ and $\lambda_2$ and the intersection with $C'$) over which the ramification is tame.

\endroster 

\vskip.6cm
 
\epsfysize=2.6cm

\centerline{\epsfbox{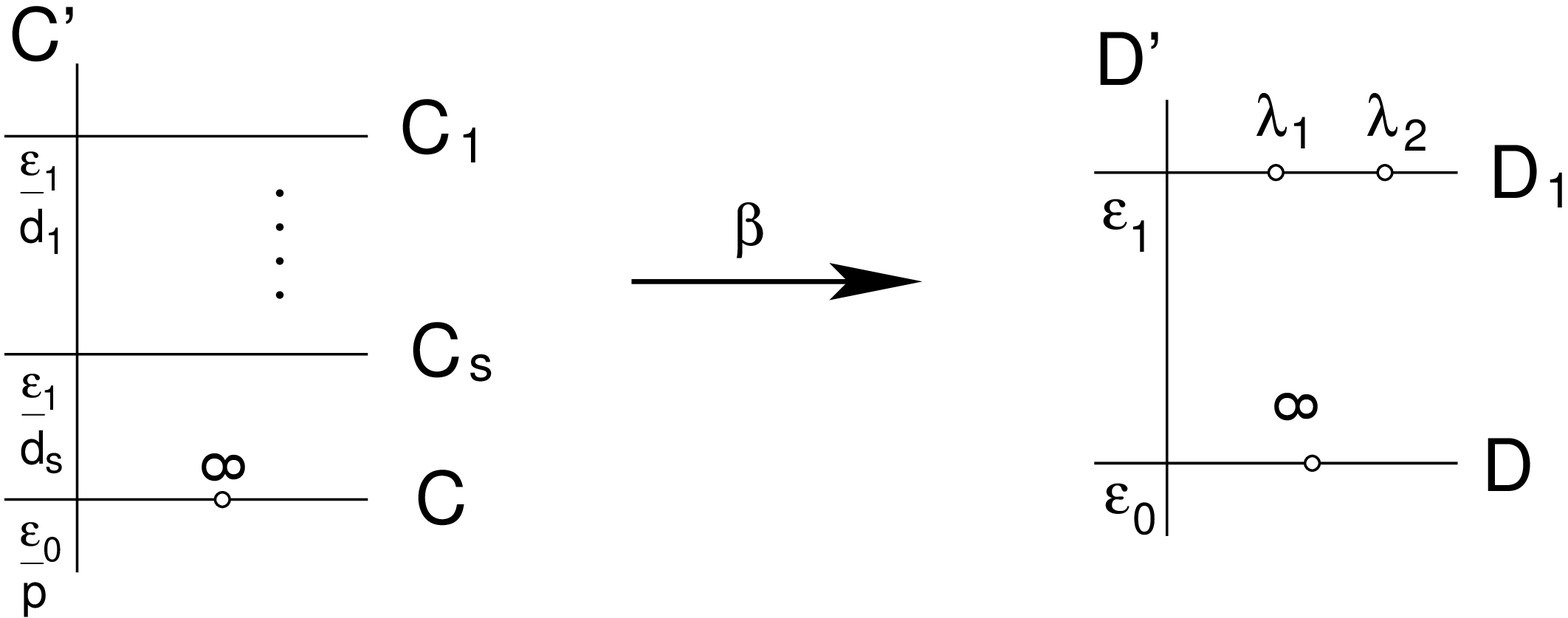}}

\vskip.6cm 
\endproclaim 

 The theorems 6.1,7.1,7.3 and  7.5 allow us to completely classify the reduction type of the  minimal semi-stable model $\Cal X@>>>\Cal Y$ (separating the ramified fibers) for a polynomial cover $\beta$ over $K$ of degree $p$, under the assumption of simple reduction of its branch locus $B$. More precisely, the above results show that the behaviour of this model  essentially depends only on the ramification data of the cover and on the thicknesses of the singular points of the special fiber of the stable (separating) model $\Cal B$ associated to the pointed curve $(\bold P^1,B)$.  

\specialhead \S8 An example\endspecialhead

We will end this paper with an explicit example. Let $K$ be a $5$-adic field, i.e. a finite extension of the field $\bold Q_5$. We will describe the semi-stable model for a polynomial cover $\beta:\bold P^1_K@>>>\bold P^1_K$ of degree $5$ having the following ramification data:
\roster
\item The cover is unramified outside the set $\{\infty,0,\lambda_1,\lambda_2\}$ with $\lambda_1,\lambda_2\in K^*$ and $\lambda_1\neq\lambda_2$. 
\item There is only one point above $\infty$. 
\item There are three points above $0$, one of them, denoted by $P_0$ has ramification index $3$, while the others are unramified.
\item There are four points above $\lambda_1$ (resp. above $\lambda_2$), one of them, denoted by $P_1$ (resp. $P_2$) has ramification index $2$ while the other three are unramified.
\endroster
Up to equivalence (cf. \S2), and since the behaviour of the ramification above $\lambda_1$ and $\lambda_2$ is the same, we can reduce to the case $\lambda_1=1$ and $\lambda_2=\lambda\in R-\{0,1\}$. The results of this paper allow us to describe the semi-stable model for $\beta$ without any direct computation. More precisely, if the branch locus has good reduction, i.e. if $v(\lambda)=v(\lambda-1)=0$ then we can apply Theorem 6.1, which leads to the following model:

\midinsert
\epsfysize=2.8cm
\centerline{\epsfbox{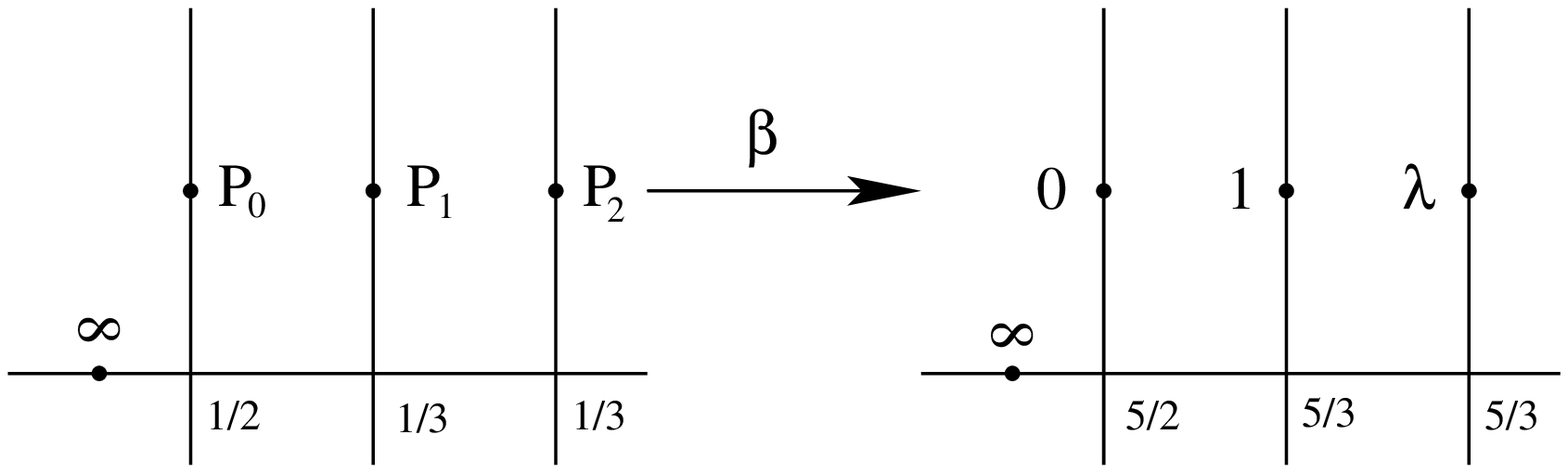}}
\endinsert

The branch locus of the cover $\beta$ has bad reduction if either $v(\lambda)>0$ or $v(\lambda-1)>0$. In the first case, there are three possibilities, leading to four different semi-stable models: if $\epsilon=v(\lambda)<5$, we can apply Theorem 7.1 and obtain the following model:

\newpage 

\midinsert
\epsfysize=2.8cm
\centerline{\epsfbox{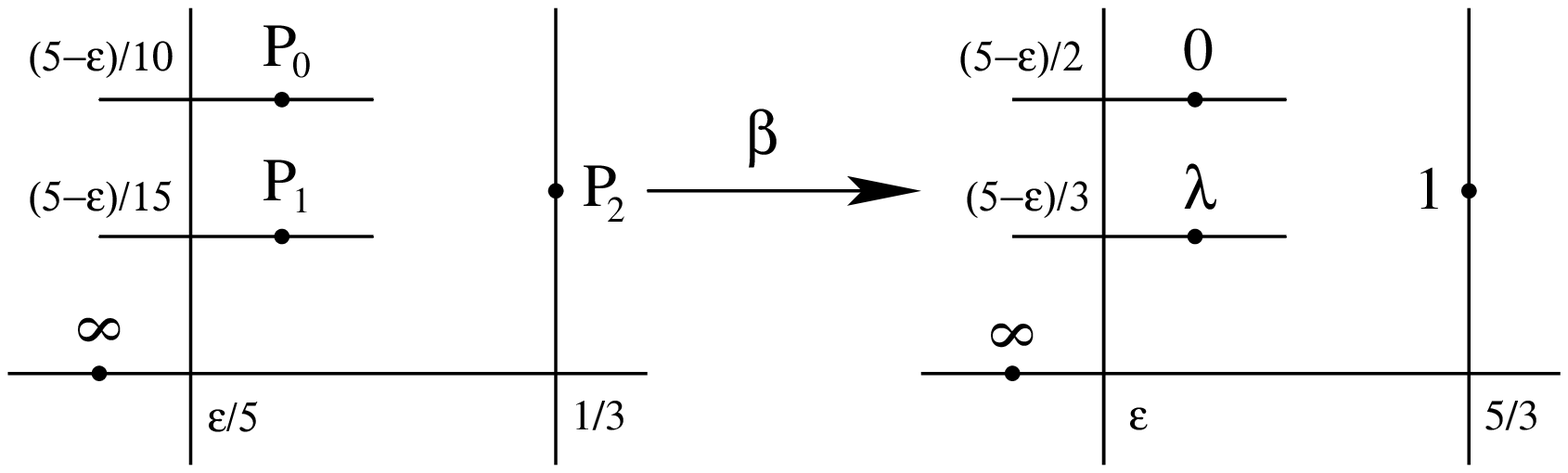}}
\endinsert

\noindent As it is shown in the following picture, and applying Teorem 7.3, the simplest case occurs for $\epsilon=5$:

\ 

\epsfysize=2.8cm
\centerline{\epsfbox{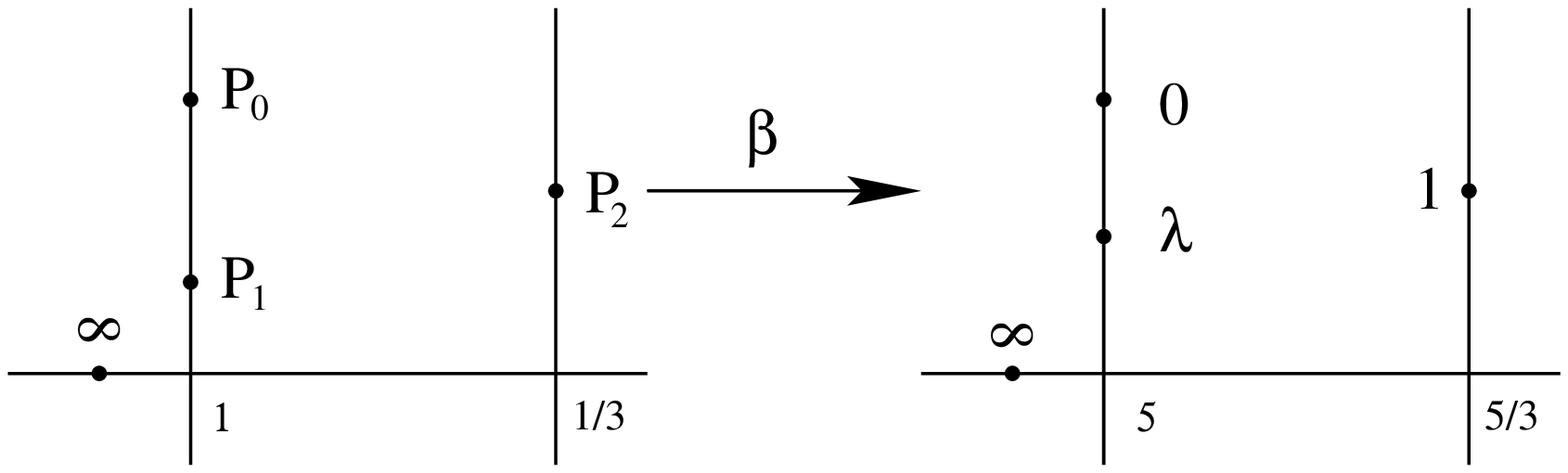}}

\ 

\noindent Finally, for $\epsilon>5$ we find two possibilities, depending on the partitions of the ramification indices $(1,1,3)$ and $(1,1,1,2)$ of the fibers above $0$ and $\lambda$ (cf. Theorem 7.5). The first partition is given by $S_0=\{\{1,3\},\{1\}\}$ and $S_1=\{\{1,1,2\},\{1\}\}$, which leads to the following semi-stable model:

\epsfysize=2.8cm
\centerline{\epsfbox{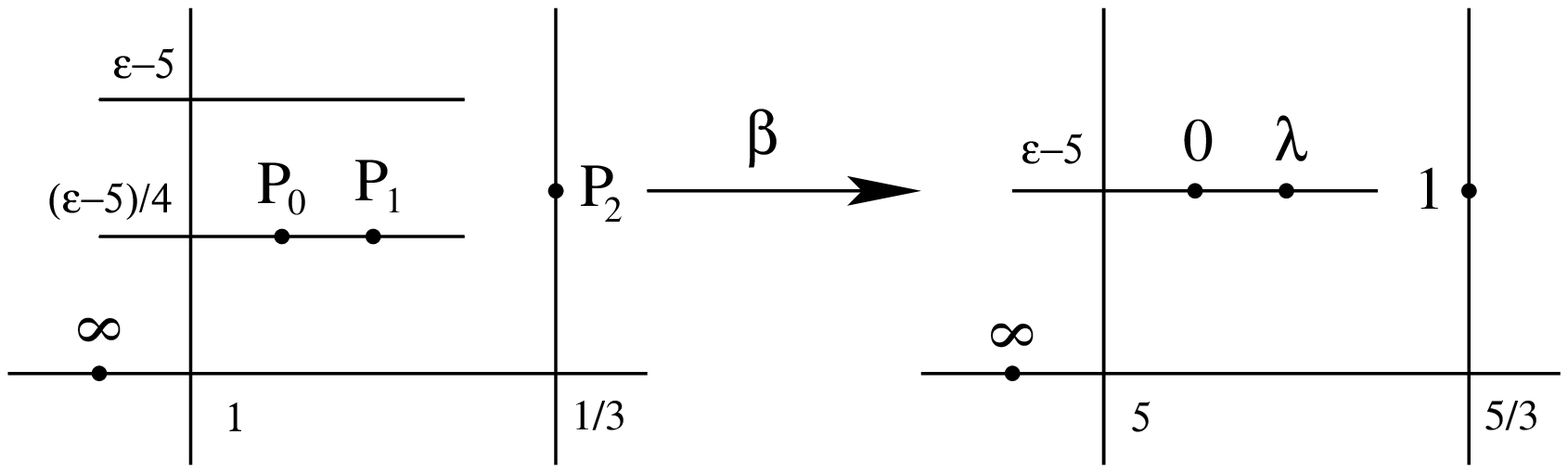}}

\noindent If we consider the second partition $S_0=\{\{3\},\{1,1\}\}$ and $S_1=\{\{1,1,1\},\{2\}\}$, we then obtain the following model:

\midinsert
\epsfysize=2.8cm
\centerline{\epsfbox{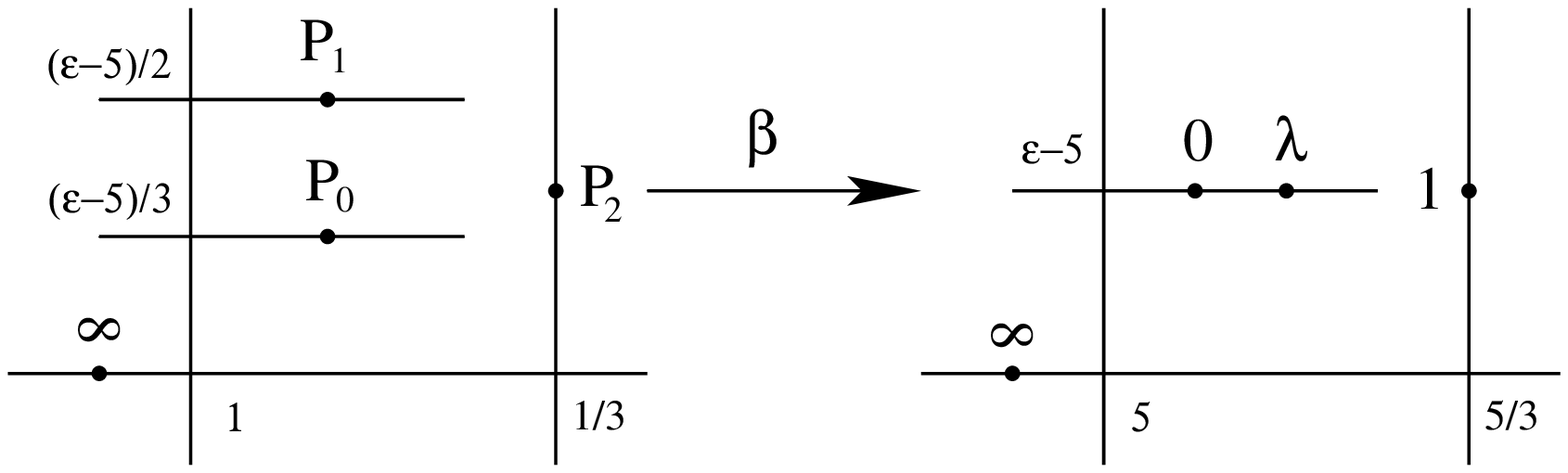}}
\endinsert

In order to completely classify the possible semi-stable models for $\beta$, we have to study the last case of bad reduction of the branch locus, i.e. when $\epsilon=v(\lambda-1)>0$. As before, there will be three different cases, depending on $\epsilon$, leading to four different situations. First of all, for $\epsilon<\frac52$ we obtain the following semi-stable model (cf. Theorem 7.1):

\epsfysize=2.8cm
\centerline{\epsfbox{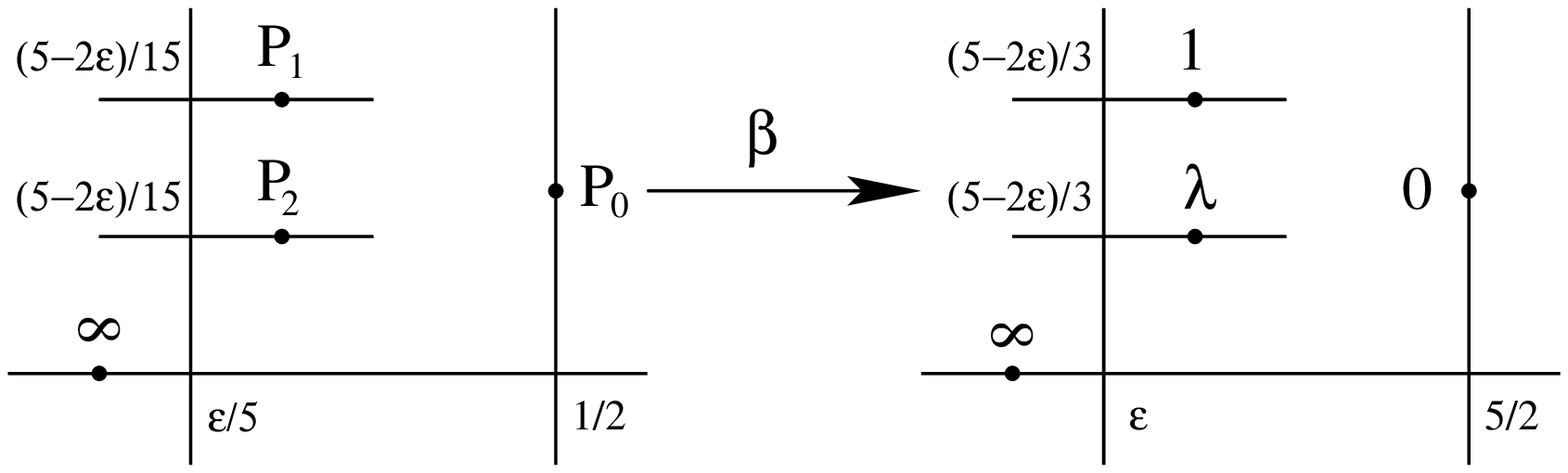}}

\noindent The next picture describes the case $\epsilon=\frac52$ (cf. Theorem 7.3):

\midinsert
\epsfysize=2.8cm
\centerline{\epsfbox{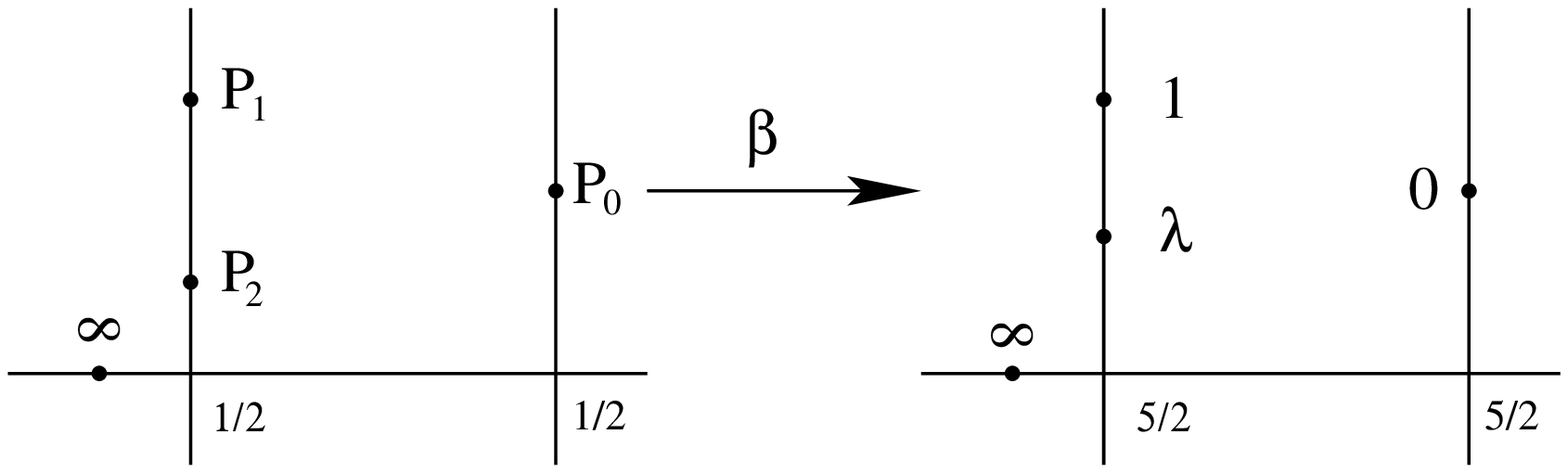}}
\endinsert

\noindent For $\epsilon>\frac52$, according to Theorem 7.5, there are two possible partitions of the ramification indices above $1$ and $\lambda$. The first is given by $S_1=S_\lambda=\{\{1,2\},\{1\},\{1\}\}$ and leads to the following semi-stable model:

\midinsert
\epsfysize=2.8cm
\centerline{\epsfbox{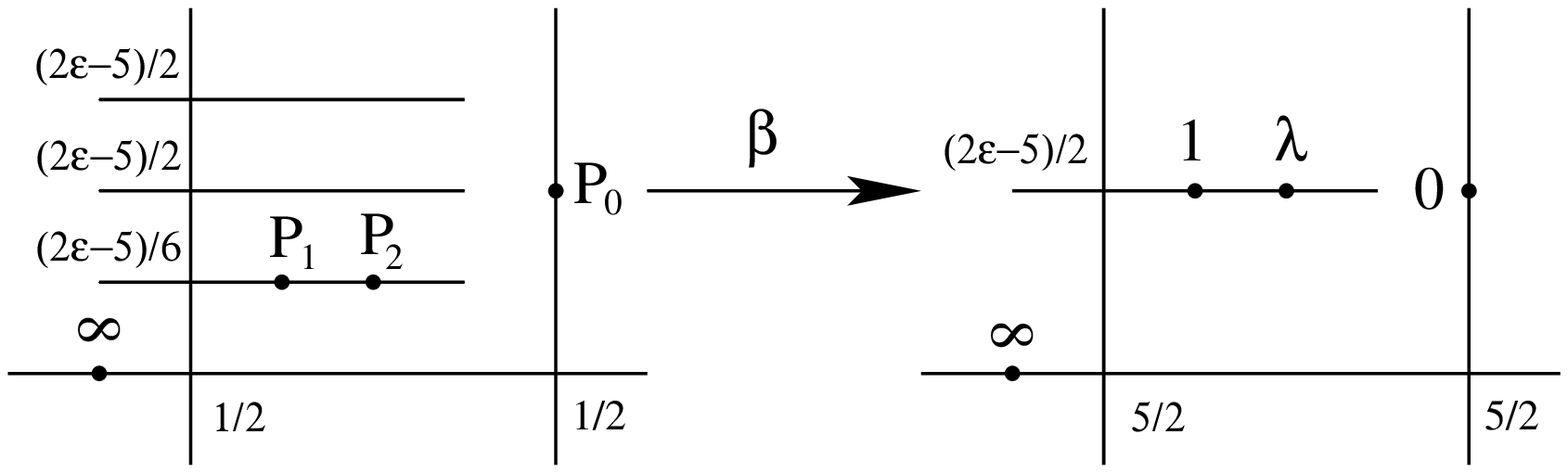}}
\endinsert

\noindent Finally, the next picture describes the semi-stable model associated to the second partition $S_1=S_\lambda=\{\{2\},\{1,1\},\{1\}\}$:

\midinsert
\epsfysize=2.8cm
\centerline{\epsfbox{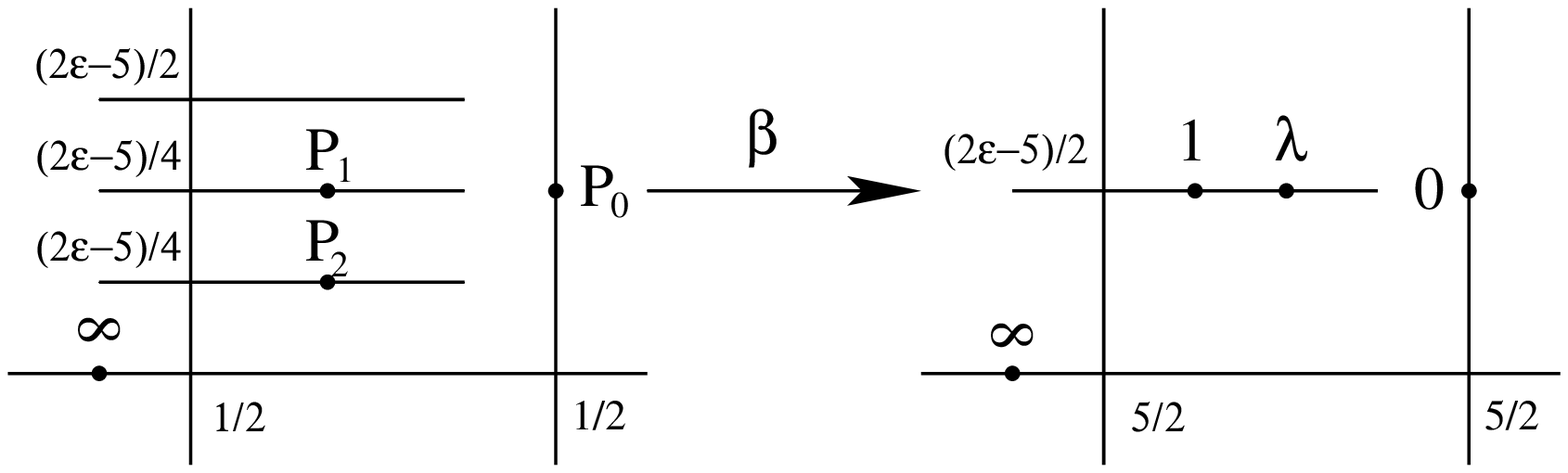}}
\endinsert


\enddocument

\end